\documentclass[12pt,a4paper,reqno]{amsart}

\usepackage{amssymb}
\usepackage{amsmath}
\usepackage{latexsym}
\usepackage{exscale}
\usepackage[latin1]{inputenc}
\usepackage{epsfig}
\usepackage{verbatim}
\usepackage{graphicx}
\usepackage{subfigure}

\newcommand{\norm}[1]{\left\Vert#1\right\Vert}
\newcommand{\abs}[1]{\left| #1 \right|}
\newcommand{\ip}[2]{\langle #1 , #2 \rangle}

\numberwithin{equation}{section} 

\newtheorem{thm}{Theorem}[section] 

\newtheorem{defi}[thm]{Definition}
\newtheorem{rem}[thm]{Remark}

\newtheorem{lem}[thm]{Lemma}

\begin{document}

\allowdisplaybreaks

\title[SHEAR ANISOTROPIC INHOMOGENEOUS SPACES]
 {\large SHEAR ANISOTROPIC INHOMOGENEOUS BESOV and TRIEBEL-LIZORKIN SPACES in $\mathbb{R}^d$}

\author{Daniel Vera}

\address{Daniel Vera
\\
Departamento de Matem\'aticas
\\
Universidad Aut\'onoma de Madrid
\\
28049 Madrid, Spain}

\email{daniel.vera@uam.es}


\date{\today}



\begin{abstract}
We define distribution spaces in $\mathbb{R}^d$ via $\ell^q(L^p)$ or
$L^p(\ell^q)$ norms of a sequence of convolutions of
$f\in\mathcal{S}'$ with smooth functions, the shearlet system. Then,
we define associated sequence spaces and prove characterizations. We
also show a reproducing identity in $\mathcal{S}'$. Finally, we
prove Sobolev-type embeddings within the shear anisotropic
inhomogeneous spaces and embeddings between (classical dyadic)
isotropic inhomogeneous spaces and shear anisotropic inhomogeneous
spaces.
\end{abstract}

\maketitle


\vskip 1cm
\section{Introduction.}\label{S:Intro}
The traditional (separable) multidimensional wavelets are built from
tensor products of $1$-dimensional wavelets. Hence, wavelets are not
very efficient in ``sensing" the geometry of lower dimension
discontinuities since the number of wavelets remains the same across
scales. In applications it may be desirable to be able to detect
more orientations having still a basis-like representation. In
recent years there have been attempts to achieve this sensitivity to
more orientations. Some of them include the directional wavelets or
filterbanks \cite{AMV}, \cite{BS}, the curvelets \cite{CaDo00} and
the contourlets \cite{DoVe}, to name just a few.

In \cite{GLLWW}, Guo, Lim, Labate, Weiss and Wilson, introduced the
wavelets with composite dilation. This type of representation takes
full advantage of the theory of affine systems on $\mathbb{R}^d$ and
therefore, unlike other representations, provides a natural
transition from the continuous representation to the discrete
(basis-like) setting, as in the case of wavelets. Based on the
wavelets of composite dilation theory, the shearlets system provides
Parseval frames for $L^2(\mathbb{R}^d)$ or subspaces of it,
depending on the discrete sampling of parameters (see \cite[Section
5.2]{GLLWW} or \cite{GKL05}).

There are at least two other ways to define smoothness spaces with
shearlets. The sophisticated theory of coorbit spaces uses auxiliary
sets and spaces of functions to define Banach spaces called
\emph{shearlet coorbit spaces}, as developed by Dahlke, Häuser,
Kutyniok, Steidl and Teschke in \cite{DHST11,DKST,DST11,DST12}. A
closer-in-spirit approach is that of Labate, Mantovani, and Negi in
\cite{LMN2012}, in which a general theory of decomposition spaces is
used on the shearlets system to define (quasi-)Banach spaces called
\emph{shearlet smoothness spaces}. Both of these approaches are
related to Besov spaces. This can easily be seen by their definition
with the shearlets coefficients in $\ell^{p,q}$ norms. One main
difference of this article with the \emph{shearlet coorbit spaces}
and the \emph{shearlet smoothness spaces} is that we define the
\textbf{shear anisotropic inhomogeneous Besov and Triebel-Lizorkin
\emph{distribution} spaces} from $\ell^q(L^p)$ and $L^p(\ell^q)$
norms, respectively, of a sequence of convolutions of
$f\in\mathcal{S}'$ with a shear anisotropic dilated functions in
$\mathcal{S}$. Then, we define the associated \emph{sequence} spaces
and prove characterization. The line of argumentation follows the
work of Frazier and Jawerth in \cite{FJ85} and \cite{FJ90}. In
\cite{DHST11,DKST,DST11,DST12} it is stated the existence of a
bounded linear reconstruction operator and it is also proved the
embedding of certain \emph{shearlet coorbit spaces} into (sums of)
homogeneous Besov spaces with $1\leq p=q\leq\infty$, for a certain
smoothness parameter. Regarding \cite{LMN2012}, the embeddings
between the \emph{shearlet smoothness spaces} and the classical
Besov spaces are in both directions for a certain, more general, set
of parameters. We prove a reproducing identity with convergence in
$\mathcal{S}'$ (previously possible only for $L^2$). We also prove
embeddings within shear anisotropic spaces and between shear
anisotropic and classical spaces.

We remind the reader the definitions of classical spaces (as appear
in \cite{FJ85} and \cite{FJ90}) and which will also be used in the
embedding results.

Let $\varphi,\Phi\in\mathcal{S}(\mathbb{R}^d)$ satisfy
\begin{equation}\label{e:supp_dyad_fcn}
    \text{supp }\hat{\varphi}\subseteq\{\xi\in\hat{\mathbb{R}}^d:\frac{1}{2}\leq\abs{\xi}\leq
    2\}, \;\;\;\;\;\;\;\;\; \text{supp }\Phi\subseteq\{\xi\in\hat{\mathbb{R}}^d:\abs{\xi}\leq 2\},
\end{equation}
and
\begin{equation}\label{e:postvty_dyad_fcn}
\abs{\varphi(\xi)}\geq c>0, \text{ if }
\frac{3}{5}\leq\abs{\xi}\leq\frac{5}{3}, \;\;\;\;\;\;\;\;\;
\abs{\Phi(\xi)}\geq c>0, \text{ if } \abs{\xi}\leq\frac{5}{3}.
\end{equation}
For $\nu\in\mathbb{N}\cup \{0\}$ and $k\in\mathbb{Z}^d$, the dyadic
dilation is $\varphi_{2^\nu I}(x):=2^{d\nu}\varphi(2^\nu x)$, where
$I$ is the identity matrix. Identifying $Q$ with $(\nu,k)$, the
dyadic dilation normalized in $L^2$ is
$\varphi_Q=\varphi_{\nu,k}:=2^{d\nu/2}\varphi(2^\nu x)$. For
$Q_0=[0,1]^d$ let $Q_{\nu,k}=2^{-\nu}(Q_0+k)$. Let $\mathcal{D}_+$
denote the set of dyadic cubes $\{Q_{\nu,k}: \nu\in\mathbb{N}_0,
k\in\mathbb{Z}^d\}$ and $\mathcal{D}^\nu_+$ denote the set of all
dyadic cubes at scale $\nu$, \emph{i.e.},
$\mathcal{D}^\nu_+=\{Q_{\nu,k}: k\in\mathbb{Z}^d\}$. Let $\chi_Q$ be
the characteristic function of $Q$ and
$\tilde{\chi}_Q=\abs{Q}^{-\frac{1}{2}}\chi_Q$ the $L^2$-normalized
characteristic function of $Q$.

For $\alpha\in\mathbb{R}$, $0<p,q\leq\infty$, the inhomogeneous
Besov distribution space $\mathbf{B}^{\alpha,q}_p$ is the set of all
$f\in\mathcal{S}'$ such that
\begin{equation}\label{e:BsvSpcs_dyad}
\norm{f}_{\mathbf{B}^{\alpha,q}_p}=\norm{\Phi\ast f}_{L^p} +
\left(\sum_{\nu=0}^\infty (2^{\nu\alpha}\norm{\varphi_{2^\nu I}\ast
f}_{L^p})^q\right)^{1/q}<\infty,
\end{equation}
and the inhomogeneous Besov sequence space $\mathbf{b}^{\alpha,q}_p$
is the set of all complex-valued sequences
$\mathbf{s}=\{s_Q\}_{Q\in\mathcal{D}_+}$ such that
\begin{equation}\label{e:BsvSeq_dyad}
\norm{\mathbf{s}}_{\mathbf{b}^{\alpha,q}_p}=
\left(\sum_{\nu=0}^\infty
\left(\sum_{Q\in\mathcal{D}^\nu_+}[\abs{Q}^{-\frac{\alpha}{d}+\frac{1}{p}-\frac{1}{2}}\abs{s_Q}]^p\right)^{q/p}\right)^{1/q}<\infty.
\end{equation}
For $\alpha\in\mathbb{R}$, $0<p<\infty$ and $0<q\leq\infty$, the
inhomogeneous Triebel-Lizorkin distribution space
$\mathbf{F}^{\alpha,q}_p$ is the set of all $f\in\mathcal{S}'$ such
that
\begin{equation}\label{e:T-LSpcs_dyad}
\norm{f}_{\mathbf{F}^{\alpha,q}_p}=\norm{\Phi\ast f}_{L^p} +
\norm{\left(\sum_{\nu=0}^\infty (2^{\nu\alpha}\abs{\varphi_{2^\nu
I}\ast f})^q\right)^{1/q}}_{L^p}<\infty,
\end{equation}
and the inhomogeneous Triebel-Lizorkin sequence space
$\mathbf{f}^{\alpha,q}_p$ is the set of all complex-valued sequences
$\mathbf{s}=\{s_Q\}_{Q\in\mathcal{D}_+}$ such that
\begin{equation}\label{e:T-LSeq_dyad}
\norm{f}_{\mathbf{f}^{\alpha,q}_p}=
\norm{\left(\sum_{Q\in\mathcal{D}_+}
(\abs{Q}^{-\frac{\alpha}{d}}\abs{s_Q}\tilde{\chi}_Q(\cdot))^q\right)^{1/q}}_{L^p}<\infty.
\end{equation}
Let $\psi$ and $\Psi$ be defined as $\varphi$ and $\Phi$,
respectively. The main results in \cite{FJ85} and \cite{FJ90} are
that the distribution spaces (\ref{e:BsvSpcs_dyad}) and
(\ref{e:T-LSpcs_dyad}) can be characterized by the corresponding
sequence spaces (\ref{e:BsvSeq_dyad}) and (\ref{e:T-LSeq_dyad}) via
the analysis and synthesis operators, formally defined as
$$S_{\Phi,\varphi}
f=\{\{\ip{f}{\Phi(\cdot-k)}\}_{k\in\mathbb{Z}^d},\{\ip{f}{\varphi_Q}\}_{Q\in\mathcal{D}_+}\},$$
$$\text{and }\;\; T_{\Psi,\psi} \mathbf{s}(\cdot)=\sum_{k\in\mathbb{Z}^d}s_k\Psi(\cdot-k)+\sum_{Q\in\mathcal{D}_+}s_Q
\psi_Q(\cdot).
$$
The characterization of the Triebel-Lizorkin spaces in \cite{FJ90}
uses the next definitions (which will be used for the embedding
results): identify $Q$ and $P$ with $(\nu,k)$ and $(\nu,k')$,
respectively. For all $r>0$, $N\in\mathbb{N}$ and $\nu\geq 0$,
define
$$(s^\ast_{r,N})_Q:=\left(\sum_{Q\in\mathcal{D}_+}\frac{\abs{s_P}^r}{(1+2^\nu\abs{x_Q-x_P})^N}\right)^{1/r},$$
where $x_Q=2^{-\nu}k$ is the ``lower left corner" of $Q_{\nu,k}$.
Define also
$\mathbf{s}^\ast_{r,N}:=\{(s^\ast_{r,N})_Q\}_{Q\in\mathcal{D}_+}$.
With an additional assumption on the functions $\psi$, $\Psi$,
$\varphi$ and $\Phi$, a reproducing identity with convergence in
$\mathcal{S}'$ is also proved in \cite[Lemma 2.1]{FJ85}. This means
that $f=T_{\Psi,\psi}\circ S_{\Phi,\phi} f$ with convergence in
$\mathcal{S}'$. We remind the reader that our notation is slightly
different from \cite{FJ85} and \cite{FJ90} with respect to the
dilations and, therefore, our notation for the analysis and
synthesis operators is different from \cite{FJ85} and \cite{FJ90}
(see (2.14) in \cite{FJ85}).

The reader can now compare (\ref{e:BsvSpcs_dyad}) and
(\ref{e:BsvSeq_dyad}) with (\ref{e:def_Bsv_Shrlts}) and
(\ref{e:def_BsvSeq_Shrlts}) for the Besov-like spaces and
(\ref{e:T-LSpcs_dyad}) and (\ref{e:T-LSeq_dyad}) with
(\ref{e:def_T-L_func_spcs}) and (\ref{e:def_T-L_seq_spcs}) for the
Triebel-Lizorkin-like spaces. See Section \ref{S:Shrlts_HD} for the
definitions of shear anisotropic dilations.

Some results for shear anisotropic inhomogeneous Triebel-Lizorkin
spaces for $d=2$ appeared first in \cite{Ve12}.

The outline of the paper is as follows. In Section \ref{S:Shrlts_HD}
we introduce the ``shearlets on the cone" system and set notation.
In Section \ref{S:Basic_Results} we give two basic lemmata regarding
almost orthogonality and present classical results slightly
modified. In Sections \ref{S:AB-Bsv} and \ref{S:AB-T-L} we define
the shear anisotropic inhomogeneous (Besov and Triebel-Lizorkin,
respectively) distributions and sequence spaces and prove
characterization. A reproducing identity in $\mathcal{S}'$ is proved
in Section \ref{S:ReprId_SmthShrlts}. In Section \ref{S:Embeddings}
we prove Sobolev-type embeddings within the shear anisotropic
inhomogeneous spaces and embeddings between (classical dyadic)
isotropic inhomogeneous spaces and shear anisotropic inhomogeneous
spaces for certain smoothness parameters. In this section we also
prove that there exist sequences of non-vanishing functions in one
of the classical or shear anisotropic spaces that vanish in the norm
of the other space. Proofs for Section \ref{S:Basic_Results} and for
some results in Section \ref{S:AB-T-L} are given in Section
\ref{S:Proofs}.

\vskip 1cm
\section{Shearlets on the cone and notation}\label{S:Shrlts_HD}
Define the cone aligned with the $\xi_1$ axis as
\begin{equation}\label{e:Cone_Domain}
\mathcal{D}^{(1)}=\{(\xi_1,\ldots,\xi_d)\in\hat{\mathbb{R}}^d:
    \abs{\xi_1}\geq\frac{1}{8}, \abs{\frac{\xi_\frak{d}}{\xi_1}}\leq 1, \frak{d}=2,\ldots,d \}.
\end{equation}
Let $\hat{\psi}_1,\hat{\psi}_2\in C^\infty(\mathbb{R})$ with
$\text{supp }\hat{\psi}_1\subset
[-\frac{1}{2},-\frac{1}{16}]\cup[\frac{1}{16},\frac{1}{2}]$ and
$\text{supp }\hat{\psi}_2\subset [-1,1]$ such that
\begin{equation}\label{e:Discrt_Shrlt_Cond_Cone_1}
\sum_{j\geq 0}\abs{\hat{\psi}_1(2^{-2j}\omega)}^2 =1, \;\;\;
\text{for }\abs{\omega}\geq \frac{1}{8}
\end{equation}
and
\begin{equation}\label{e:Discrt_Shrlt_Cond_Cone_2}
\abs{\hat{\psi}_2(\omega-1)}^2+\abs{\hat{\psi}_2(\omega)}^2+\abs{\hat{\psi}_2(\omega+1)}^2=1,
\;\;\; \text{for } \abs{\omega}\leq 1.
\end{equation}
It follows from (\ref{e:Discrt_Shrlt_Cond_Cone_2}) that, for $j\geq
0$,
\begin{equation}\label{e:Discrt_Shrlt_Cond_Cone_3}
\sum_{\ell=-2^j}^{2^j} \abs{\hat{\psi}_2(2^j\omega-\ell)}^2=1,
\;\;\; \text{for }\abs{\omega}\leq 1.
\end{equation}
For a scale $j\geq 0$ the anisotropic dilation matrices are defined
as
$$A^j_{(1)}= \left(
\begin{array}{cccc}
  4^j    & 0      & \ldots & 0\\
  0      & 2^j    & \ldots & 0\\
  \vdots & \vdots & \vdots & \vdots\\
  0      & 0      & \ldots & 2^j\\
\end{array}
\right), \;\;\; \ldots \;\;\; , A^j_{(d)}=\left(
\begin{array}{cccc}
  2^j    & 0      & \ldots & 0\\
  0      & 2^j      & \ldots & 0\\
  \vdots & \vdots & \vdots & \vdots\\
  0      & 0      & \ldots & 4^j\\
\end{array}
\right),
$$
and for $\ell=(\ell_1,\ldots,\ell_{d-1})$ with
$-2^j\leq\ell_i\leq2^j$, $i=1,\ldots,d-1$, the $d\times d$ shear
matrices are defined as
$$B^{[\ell]}_{(1)}= \left(
\begin{array}{cccc}
  1      & \ell_1 & \ldots & \ell_{d-1}\\
  0      & 1      & \ldots & 0\\
  \vdots & \vdots & \vdots & \vdots\\
  0      & 0      & \ldots & 1\\
\end{array}
\right), \;\;\; \ldots, B^{[\ell]}_{(d)}=\left(
\begin{array}{cccc}
  1      & 0      & \ldots & 0\\
  0      & 1      & \ldots & 0\\
  \vdots & \vdots & \vdots & \vdots\\
\ell_1   & \ell_2 & \ldots & 1\\
\end{array}
\right).
$$
To shorten notation we will write $\abs{[\ell]}\preceq 2^j$ instead
of $\abs{\ell_i}\leq 2^j$, $i=1,\ldots,d-1$, and $\abs{[\ell]}=2^j$
when $\abs{\ell_i}=2^j$ for at least one $i=1,2,\ldots,d-1$. Define
$\hat{\psi}^{(1)}(\xi):=\hat{\psi}_1(\xi_1)\prod_{\frak{d}=2}^{d}\hat{\psi}_2(\frac{\xi_\frak{d}}{\xi_1})$.
Since $\xi A_{(1)}^{-j} B_{(1)}^{[-\ell]}=(4^{-j}\xi_1,
-4^{-j}\xi_1\ell_1+2^{-j}\xi_2, \ldots,
-4^{-j}\xi_1\ell_{d-1}+2^{-j}\xi_d)$, from
(\ref{e:Discrt_Shrlt_Cond_Cone_1}) and
(\ref{e:Discrt_Shrlt_Cond_Cone_3}) it follows that
\begin{eqnarray}\label{e:PrsvlFrm_Prop_Shrlts}
\nonumber
    & &\!\!\!\!\!\!\!\!\!\!\!\!\!\!\!\!\!\!\!\!\!\!\!\!\!\!\!\!\!\!\!\!\!\!\!\!\!\!\!\!\!\!\!\!
        \sum_{j\geq 0}\sum_{\abs{[\ell]}\preceq2^j}
        \abs{\hat{\psi}^{(1)}(\xi A^{-j}_{(1)}B^{[-\ell]}_{(1)})}^2 \\
\nonumber
    &=& \sum_{j\geq0}\sum_{\abs{\ell_1},\ldots,\abs{\ell_{d-1}}\leq2^j}
        \abs{\hat{\psi}_1(2^{-2j}\xi_1)}^2\prod_{\frak{d}=2}^d\abs{\hat{\psi}_2(2^j\frac{\xi_\frak{d}}{\xi_1}-\ell_{\frak{d}-1})}^2 \\
\nonumber
    &=& \sum_{j\geq0} \abs{\hat{\psi}_1(2^{-2j}\xi_1)}^2\sum_{\abs{\ell_1},\ldots,\abs{\ell_{d-2}}\leq2^j}
        \prod_{\frak{d}=2}^{d-1}\abs{\hat{\psi}_2(2^j\frac{\xi_\frak{d}}{\xi_1}-\ell_{\frak{d}-1})}^2\\
\nonumber
    &\vdots&  \\
    &=& 1,
\end{eqnarray}
for $\xi=(\xi_1,\ldots,\xi_d)\in\mathcal{D}^{(1)}$ and which we will
call the \textbf{Parseval frame condition} (for the cone
$\mathcal{D}^{(1)}$). Since $\text{supp }\hat{\psi}^{(1)}\subset
[-\frac{1}{2},\frac{1}{2}]^d$, (\ref{e:PrsvlFrm_Prop_Shrlts})
implies that the shearlet system
\begin{equation}\label{e:Shrlt_Sys_Cone}
\{\psi_{j,\ell,k}^{(1)}(x)= \abs{\text{det
}A_{(1)}}^{j/2}\psi^{(1)}(B^{[\ell]}_{(1)} A^j_{(1)} x-k): j\geq 0,
\abs{[\ell]}\preceq2^j, k\in\mathbb{Z}^d\},
\end{equation}
is a Parseval frame for $L^2((\mathcal{D}^{(1)})^\vee)=\{f\in
L^2(\mathbb{R}^d): \text{supp }\hat{f}\subset\mathcal{D}^{(1)}\}$
(see \cite{GLLWW}, Subsection 5.2.1). This means that
$$\sum_{j\geq 0}\sum_{\abs{[\ell]}\preceq 2^j}\sum_{k\in\mathbb{Z}^d}
    \abs{\ip{f}{\psi_{j,\ell,k}^{(1)}}}^2 = \norm{f}^2_{L^2(\mathbb{R}^d)},$$
for all $f\in L^2(\mathbb{R}^d)$ such that $\text{supp
}\hat{f}\subset \mathcal{D}^{(1)}$. 
There are several examples of functions $\psi_1,\psi_2$ satisfying
the properties described above (see \cite{GL07}). Since
$\hat{\psi}^{(1)}\in C^\infty_c(\hat{\mathbb{R}}^d)$,  there exists
$C_N$ such that $\abs{\psi^{(1)}(x)}\leq C_N (1+\abs{x})^{-N}$ for
all $N\in\mathbb{N}$. The geometric properties of the shearlets
system in $\mathcal{D}^{(1)}$ are more evident by observing that
$$\text{supp }(\psi^{(1)}_{j,\ell,k})^\wedge\subset
\{\xi\in\hat{\mathbb{R}}^d : \abs{\xi_1}\in [2^{2j-4},2^{2j-1}],
\abs{2^j\frac{\xi_\frak{d}}{\xi_1}-\ell_{\frak{d}-1}}\leq 1,
\frak{d}=2,\ldots,d\}.$$

\noindent One can also construct a shearlets system for any cone
\begin{equation*}
\mathcal{D}^{(i)}=\{\xi=(\xi_1,\ldots,\xi_i,\ldots,\xi_d)\in\hat{\mathbb{R}}^d:
    \abs{\xi_i}\geq\frac{1}{8}, \abs{\frac{\xi_\frak{d}}{\xi_i}}\leq 1, \frak{d}\not=i\},
\end{equation*}
by defining
$\hat{\psi}^{(i)}(\xi)=\hat{\psi}_1(\xi_i)\prod_{\frak{d}\not=i}\hat{\psi}_2(\frac{\xi_\frak{d}}{\xi_i})$
and choosing correspondingly the anisotropic and shear matrices
$A_{(i)}^j$ and $B_{(i)}^{[\ell]}$.

Let $\hat{\Psi}\in C^\infty_c(\mathbb{R}^d)$, with $\text{supp
}\hat{\Psi}\subset [-\frac{1}{4},\frac{1}{4}]^d$ and
$\abs{\hat{\Psi}}=1$ for $\xi\in
[-\frac{1}{8},\frac{1}{8}]^d=\mathcal{R}$, be such that
\begin{equation}\label{e:ShrltPrsvlFrm}
   \abs{\hat{\Psi}(\xi)}^2\chi_\mathcal{R}(\xi)
        + \sum_{\frak{d}=1}^d\sum_{j\geq 0}\sum_{\abs{[\ell]}\preceq2^j}\abs{\hat{\psi}^{(\frak{d})}(\xi A^{-j}_{(\frak{d})}B^{[-\ell]}_{(\frak{d})})}^2
        \chi_{\mathcal{D}^{(\mathfrak{d})}}(\xi)=1,
\end{equation}
for all $\xi\in\hat{\mathbb{R}}^d$. This implies that one can
construct a Parseval frame for $L^2(\mathbb{R}^d)$, see
\cite[Theorem 9]{LaWe09}.

Since $\mathcal{D}^{(i)}$ are orthogonal rotations of
$\mathcal{D}^{(1)}$ we will develop our results only for direction
$1$. We will incorporate the directions only in the definitions of
the spaces. Then, dropping the subindices in the matrices and the
superindices in the shearlet functions we have
\begin{equation}\label{e:BAx}
B^{[\ell]} A^j x= \left(
\begin{array}{c}
  2^{2j}x_1+2^j\ell_1x_2+\ldots+2^j\ell_{d-1}x_d      \\
  2^jx_2      \\
  \vdots \\
  2^jx_d   \\
\end{array}
\right),
\end{equation}
for every $j\geq 0$ and $\abs{[\ell]}\preceq 2^j$ and
$\psi_{j,\ell,k}(x)= \abs{\text{det }A}^{j/2}\psi(B^{[\ell]} A^j
x-k).$

\vskip 0.5cm
\subsection{Notation}\label{sS:Notation} We denote for a
matrix $M\in GL_d(\mathbb{R})$ the anisotropic dilation
$\psi_M(x)=\abs{\text{det} M}^{-1}\psi(M^{-1}x)$. We also denote
$\tilde{\psi}(x)=\overline{\psi(-x)}$. For $Q_0=[0,1)^d$, write
\begin{equation}\label{e:def_Qjlk}
Q_{j,\ell,k}=A^{-j}B^{-[\ell]}(Q_0+k),
\end{equation}
with $j\geq 0$, $\abs{[\ell]}\preceq 2^j$ and $k\in\mathbb{Z}^d$.
Therefore, $\int
\chi_{Q_{j,\ell,k}}=\abs{Q_{j,\ell,k}}=2^{-(d+1)j}=\abs{\text{det }
A}^{-j}$. We also write $\tilde{\chi}_Q(x)=\abs{Q}^{-1/2}\chi_Q(x)$.
Let $\mathcal{Q}_{AB}:=\{Q_{j,\ell,k}: j\geq 0,\abs{[\ell]}\preceq
2^j, k\in\mathbb{Z}^d\}$ and $\mathcal{Q}^{j,\ell}:=
\{Q_{j,\ell,k}:k\in\mathbb{Z}^d\}$. Then, $\mathcal{Q}^{j,\ell}$ is
a partition of $\mathbb{R}^d$. To shorten notation and clear
exposition, we will identify the multi indices $(j,\ell,k)$ and
$(i,m,n)$ with $P$ and $Q$, respectively. This way we write
$\psi_P=\psi_{j,\ell,k}$ or $\psi_Q=\psi_{i,m,n}$. Also, let $x_P$
and $x_Q$ be the ``lower left corners" $A^{-j}B^{-[\ell]}k$ and
$A^{-i}B^{-[m]}n$ of the ``cubes" $P=Q_{j,\ell,k}$ and
$Q=Q_{i,m,n}$, respectively. Let $B_r(x)$ be the Euclidean ball
centered in $x$ with radius $r$.

The elements of the shearlets system
$$\{\psi_{j,\ell,k}(x)= \abs{\text{det }A}^{j/2}\psi(B^{[\ell]} A^j x-k):
    j\geq 0, \abs{[\ell]}\preceq 2^j, k\in\mathbb{Z}^d\},$$
have Fourier transform
$$(\psi_{j,\ell,k})^\wedge(\xi)= \abs{\text{det } A}^{-j/2} \hat{\psi}(\xi A^{-j}B^{-[\ell]})
    \mathbf{e}^{-2\pi i \xi A^{-j}B^{[-\ell]} k}.$$
Using the anisotropic dilation it is also easy to verify that
$$\psi_{A^{-j}B^{-[\ell]}}(x-A^{-j}B^{-[\ell]}k)=\abs{\text{det }A}^{j/2}\psi_{j,\ell,k}(x)
    =\abs{P}^{-1/2}\psi_P(x),$$
and thus
$$\left(\psi_{A^{-j}B^{-[\ell]}}(\cdot - A^{-j}B^{-[\ell]}k)\right)^\wedge(\xi)
    =\hat{\psi}(\xi A^{-j}B^{-[\ell]})\mathbf{e}^{-2\pi i\xi A^{-j}B^{-[\ell]}k}.$$
We also have
\begin{eqnarray}\label{e:InnrProd_equiv_Conv_Shrlts}
\nonumber
  \ip{f}{\psi_P}
    &=& \ip{f}{\psi_{j,\ell,k}} \\
\nonumber
    &=& \int_{\mathbb{R}^d} f(x) \overline{\abs{\text{det }A}^{-j/2}\psi_{A^{-j}B^{[-\ell]}}(x-A^{-j}B^{[-\ell]}k)} dx \\
    &=& \abs{P}^{1/2}(f\ast\tilde{\psi}_{A^{-j}B^{[-\ell]}})(x_P).
\end{eqnarray}

We formally define the \textbf{analysis} and \textbf{synthesis
operators} as
\begin{equation}\label{e:AnlsOpShrlt}
  S_{\Psi,\psi} f = \{ \{\ip{f}{\Psi(\cdot - k)}\}_{k\in\mathbb{Z}^d} , \{\ip{f}{\psi_Q}\}_{Q\in\mathcal{Q}_{AB}} \}
\end{equation}
and
\begin{equation}\label{e:SnhtsOpShrlt}
   T_{\Psi,\psi} \mathbf{s}=\sum_{k\in\mathbb{Z}^d}s_k \Psi(\cdot -k)+\sum_{Q\in\mathcal{Q}_{AB}}s_Q \psi_Q,
\end{equation}
respectively. We remind the reader that the function related to the
low frequencies $\Psi$ (see (\ref{e:ShrltPrsvlFrm})) and associated
sequence $\{s_k\}_{k\in\mathbb{Z}^d}$ are not studied in this work
since they are already treated in the literature (see \emph{e.g.}
\cite[Section 7]{FJ85} or \cite[Section 12]{FJ90}).

\vskip 1cm
\section{Basic results}\label{S:Basic_Results}

\subsection{Almost orthogonality}\label{sS:AlmostOrthg_Basic_Results}
The next two ``almost orthogonality" results are in the form of
convolutions in the time domain. The first is between two functions
with shear anisotropic dilations and is used in both
characterizations. The second is between a function with shear
anisotropic dilation and other function with dyadic dilation and is
used in the embeddings. The last ``almost orthogonality" result is
in the Fourier domain and is used in both characterizations. The
three of these results are proved in Section
\ref{sS:Proofs_Basic_Results}.

\begin{lem}\label{l:c:l:conv_shrlts}
Let $g,h\in \mathcal{S}$. For $i=j-1,\;j,\;j+1\geq 0$, let $Q$ be
identified with $(i,m,n)$. Then, for every $N>d$, there exists a
$C_N>0$ such that
$$\abs{g_{A^{-j}B^{-\ell}}\ast h_Q(x)}\leq
\frac{C_N \abs{Q}^{-\frac{1}{2}}}{(1+2^i\abs{x-x_Q})^{N}},$$ for all
$x\in\mathbb{R}^d$.
\end{lem}

\begin{lem}\label{l:conv_shrlts-wvlts}
Let $\psi,\varphi\in\mathcal{S}$. For $j\geq 0$,
$\abs{[\ell]}\preceq 2^j$ and $k\in\mathbb{Z}^d$,
$$\int_{\mathbb{R}^d} \abs{\psi(B^{[\ell]} A^j(x-y))}\abs{\varphi(2^{2j}y)} dy
    \leq \frac{C_N\abs{P_j}}{(1+2^j\abs{x})^N},$$
for all $N>d$.
\end{lem}

\begin{lem}\label{l:Overlap_bnd}
Let $(\psi_{j,\ell,k})^\wedge$ be as in the introduction of this
section. Then, the support of $(\psi_{j,\ell,k})^\wedge$ overlaps
with the support of at most $2^{(d-1)} + 3^{(d-1)} + 6^{(d-1)}$
other shearlets $(\psi_{i,m,n})^\wedge$ for $(j,\ell)\not =(i,m)$
and all $k,n\in \mathbb{Z}^d$.
\end{lem}

\begin{rem}\label{r:l:Overlap_bnd1}
Since the translation parameters $k$ and $n$ do not affect the
support in the frequency domain, then for
$$f=T_\psi\mathbf{s}=\sum_{Q\in \mathcal{Q}_{AB}}s_Q\psi_Q
=\sum_{i\geq 0}\sum_{\abs{[m]}\preceq
2^i}\sum_{n\in\mathbb{Z}^d}s_{i,m,n}\psi_{i,m,n},$$ we formally have
that
$$(\tilde{\psi}_{A^{-j}B^{[-\ell]}}\ast f)(x)=\sum_{i=j-1}^{j+1}\sum_{m(\ell,i)}\sum_{Q\in\mathcal{Q}^{i,m}}
     s_Q(\tilde{\psi}_{A^{-j}B^{[-\ell]}}\ast \psi_Q)(x),$$
where $m(\ell,i)$ are the shear indices of those shearlets in the
Fourier domain ``surrounding" the support of
$(\tilde{\psi}_{A^{-j}B^{[-\ell]}})^\wedge$ and the sum
$\sum_{i=j-1}^{j+1}\sum_{m(\ell,i)}$ has at most $3^{(d-1)} +
3^{(d-1)} + 6^{(d-1)} +1$ terms for all $(j,\ell)$ by Lemma
\ref{l:Overlap_bnd}.
\end{rem}

\begin{rem}\label{r:l:Overlap_bnd2}
From Lemma \ref{l:Overlap_bnd} the number of shearlets in other
cones overlapping (even with shearlets within different cones) on
the Fourier domain is bounded for all scales, since their respective
systems are orthonormal rotations of the system for $\mathcal{D}^1$.
Therefore, removing the characteristic functions
$\chi_{\mathcal{D}^{(\mathfrak{d})}}$ and $\chi_\mathcal{R}$ in
(\ref{e:ShrltPrsvlFrm}) affects only the Parseval condition on the
frame.
\end{rem}

\begin{rem}\label{r:AlmstOrth_OthrDefs}
The notion of almost orthogonality has been used in the context of
shearlets in \cite{GL08} to bound the magnitude of the inner product
of more general \emph{shearlet molecules} using a dyadic parabolic
pseudo-distance. One consequence is that more general frames can be
defined by non band-limited shearlet-like functions as in
\cite{LMN2012}. In this article we only use the Euclidean distance
and prove a reproducing identity with the ``smooth Parseval frames
of shearlets".
\end{rem}

\vskip 0.5cm
\subsection{Shear anisotropic dilations and classical results}\label{sS:ShrAnstrpcClassResults_Basic_Results}
The next two Lemmata will be used to prove the characterization of
the shear anisotropic inhomogeneous Besov spaces. Their proof are in
Section \ref{sS:Proofs_Basic_Results} and are just variations of
those found in \cite{FJW}, we include them for completeness.

\begin{lem}\label{l:Smplng_BndLim_Shrlt}\textbf{[Sampling lemma]}
Let $g\in\mathcal{S}'$ and $h\in\mathcal{S}$ be such that
$$\text{\emph{supp} }\hat{g}, \text{\emph{ supp} }\hat{h}\subset [-1/2,1/2]^dB^{[\ell]} A^j,
    \;\;\;\;\;\;  j\geq 0, \abs{[\ell]}\preceq 2^j.$$
Then,
$$g\ast h= \sum_{k\in\mathbb{Z}^d} \abs{\text{\emph{det} }A}^{-j} g(A^{-j}B^{[-\ell]}k) h(x-A^{-j}B^{[-\ell]}k),$$
with convergence in $\mathcal{S}'$.
\end{lem}

\begin{lem}\label{l:Plnchl-Polya_estimate_shrlts}
\textbf{[Plancherel-Pólya]} Let $0<p\leq\infty$ and $j\geq 0$.
Suppose $g\in\mathcal{S}'$ and $\text{\emph{supp} } \hat{g}\subseteq
[-\frac{1}{2},\frac{1}{2}]^dB^{[\ell]} A^j$. Then,
$$\left(\sum_{Q\in\mathcal{Q}^{j,\ell}}\sup_{z\in
Q}\abs{g(z)}^p\right)^{1/p} \leq
C_p\abs{Q_j}^{-\frac{d}{p(d+1)}}\norm{g}_{L^p}.$$
\end{lem}

\vskip 0.5cm The next definition and result are well known and will
be used to characterize the shear anisotropic inhomogeneous
Triebel-Lizorkin distribution spaces.
\begin{defi}\label{d:H-L_Max_Fcn}
The \textbf{Hardy-Littlewood maximal function}, $\mathcal{M}f(x)$,
is given by
$$\mathcal{M}f(x)=\sup_{r>0}\frac{1}{\abs{B_r(x)}}\int_{B_r(x)} \abs{f(y)} dy,$$
for a locally integrable function $f$ on $\mathbb{R}^d$ and where
$B_r(x)$ is the ball with center in $x$ and radius $r$.
\end{defi}
\noindent It is well known that $\mathcal{M}$ is bounded on $L^p$,
$1<p\leq \infty$. It is also true that the next vector-valued
inequality holds (see \cite{FeSt71}).
\begin{thm}\label{t:Feff-Stn_Ineq} \textbf{[Fefferman-Stein]}
For $1<p<\infty$ and $1<q\leq \infty$, there exist a constant
$C_{p,q}$ such that
$$\norm{\left\{\sum_{i=1}^\infty (\mathcal{M}f_i)^q\right\}^{1/q}}_{L^p}
\leq C_{p,q}\norm{\left\{\sum_{i=1}^\infty
f_i^q\right\}^{1/q}}_{L^p},$$ for any sequence $\{f_i:
i=1,2,\ldots\}$ of locally integrable functions.
\end{thm}

\vskip 1cm
\section{Shear anisotropic inhomogeneous Besov spaces}\label{S:AB-Bsv}
After defining the spaces we will leave aside the analysis of the
low frequencies and the analysis for all directions. We follow
\cite{FJ85}.

\begin{defi}\label{d:BsvSpcs_Shrlts}
Let $\hat{\Psi}$ and $\hat{\psi}$ be such that the Parseval frame
condition (\ref{e:ShrltPrsvlFrm}) is satisfied. For $\alpha\in
\mathbb{R}$, $0<p,q\leq\infty$, the \textbf{shear anisotropic
inhomogeneous Besov \emph{distribution} space} is defined as the
collection of all $f\in\mathcal{S}'$ such that
\begin{equation}\label{e:def_Bsv_Shrlts}
\norm{f}_{\mathbf{B}^{\alpha,q}_p(AB)}:= \norm{f\ast\Psi}_{L^p}
+\left(\sum_{\mathfrak{d}=1}^d\sum_{j\geq
0}\sum_{\abs{[\ell]}\preceq 2^j}
[\abs{Q_j}^{-\alpha}\norm{f\ast\psi^{(\mathfrak{d})}_{A^{-j}_{(\mathfrak{d})}B^{-\ell}_{(\mathfrak{d})}}}_{L^p}]^q\right)^{1/q}<\infty.
\end{equation}
\end{defi}

\begin{defi}\label{d:BsvSeqSpcs_Shrlts}
For $\alpha\in \mathbb{R}$, $0<p,q\leq\infty$, the \textbf{shear
anisotropic inhomogeneous Besov \emph{sequence} space} is defined as
the collection of all complex-valued sequences
$\mathbf{s}=\{s_Q\}_{Q\in\mathcal{Q}_{AB}}$ such that
\begin{equation}\label{e:def_BsvSeq_Shrlts}
\norm{\mathbf{s}}_{\mathbf{b}^{\alpha,q}_p(AB)} :=
\left(\sum_{k\in\mathbb{Z}^d} \abs{s_k}^p\right)^{1/p} +
\left(\sum_{\mathfrak{d}=1}^d\sum_{j\geq 0}\sum_{\abs{[\ell]}\preceq
2^j}
\left(\sum_{Q\in\mathcal{Q}^{j,\ell}}[\abs{Q}^{-\alpha+\frac{d}{p(d+1)}-\frac{1}{2}}
    \abs{s_Q}]^p\right)^{q/p}\right)^{1/q}<\infty.
\end{equation}
\end{defi}

\vskip 0.5cm
\subsection{The characterization}\label{sS:Chrctzn_BsvAB_Shrlts}
We prove the boundedness of the analysis and synthesis operators
(\ref{e:AnlsOpShrlt}) and (\ref{e:SnhtsOpShrlt}) on the spaces in
Definitions \ref{d:BsvSpcs_Shrlts} and \ref{d:BsvSeqSpcs_Shrlts}.

\begin{thm}\label{t:Bndnss_T-S_Ops}
Let $\alpha\in\mathbb{R}$ and $0<p,q\leq\infty$. Then, the operators
$S_\psi : \mathbf{B}^{\alpha,q}_p(AB)\rightarrow
\mathbf{b}^{\alpha,q}_p(AB)$ and $T_\psi :
\mathbf{b}^{\alpha,q}_p(AB)\rightarrow \mathbf{B}^{\alpha,q}_p(AB)$
are well defined and bounded.
\end{thm}
\textbf{Proof.} We prove only the case $p,q<\infty$. To prove the
boundedness of $S_\psi$ assume $f\in\mathbf{B}^{\alpha,q}_p(AB)$.
From Lemma \ref{l:Smplng_BndLim_Shrlt} we have
\begin{eqnarray*}
    & & \!\!\!\!\!\!\!\!\!\!\!\!\!\!\!\!\!\!\!\!\!\!\!\!
        f\ast\tilde{\psi}_{A^{-j}B^{-\ell}}\ast\psi_{A^{-j}B^{-\ell}}(x)\\
    &=& \sum_{k\in\mathbb{Z}^d}\abs{\text{det }A}^{-j} f
        \ast\tilde{\psi}_{A^{-j}B^{-\ell}}(A^{-j}B^{-\ell}k)
        \cdot\psi_{A^{-j}B^{-\ell}}(x-A^{-j}B^{-\ell}k) \\
    &=& \sum_{k\in\mathbb{Z}^d}\abs{\text{det }A}^{-\frac{j}{2}} f
        \ast\tilde{\psi}_{A^{-j}B^{-\ell}}(A^{-j}B^{-\ell}k)
        \cdot\psi_{j,\ell,k}(x),
\end{eqnarray*}
with convergence in $\mathcal{S}'$. Identify $Q$ with $(j,\ell,k)$.
Then,
$s_Q=\ip{f}{\psi_{j,\ell,k}}=\abs{Q}^{\frac{1}{2}}f\ast\tilde{\psi}_{A^{-j}B^{-\ell}}(A^{-j}B^{-\ell}k)$
(see (\ref{e:InnrProd_equiv_Conv_Shrlts})) and Lemma
\ref{l:Plnchl-Polya_estimate_shrlts} yields
\begin{eqnarray*}
    & & \!\!\!\!\!\!\!\!\!\!\!\!\!\!\!\!\!\!\!\!\!\!\!\!
        \left(\sum_{Q\in \mathcal{Q}^{j,\ell}}
        [\abs{Q}^{-\alpha+\frac{d}{p(d+1)}-\frac{1}{2}}\abs{Q}^{\frac{1}{2}}
        \abs{f\ast\tilde{\psi}_{A^{-j}B^{-\ell}}(A^{-j}B^{-\ell}k)}]^p \right)^{1/p} \\
    &\leq& C_p \abs{Q}^{-\frac{d}{p(d+1)}} \norm{\abs{Q}^{-\alpha+\frac{d}{p(d+1)}}
        f\ast\tilde{\psi}_{A^{-j}B^{-\ell}}}_{L^p}.
\end{eqnarray*}
Therefore,
\begin{eqnarray*}
  \norm{\mathbf{s}}_{\mathbf{b}^{\alpha,q}_p(AB)}
    &\leq& C_p \left(\sum_\mathfrak{d}\sum_{j\geq 0}\sum_{\abs{[\ell]}\preceq 2^j}
        [\abs{Q}^{-\alpha}\norm{f\ast\tilde{\psi}_{A^{-j}B^{-\ell}}}_{L^p}]^q\right)^{1/q} \\
    &=& C_p \norm{f}_{\mathbf{B}^{\alpha,q}_p(AB)}.
\end{eqnarray*}

To prove the boundedness of $T_\psi$ assume
$\mathbf{s}\in\mathbf{b}^{\alpha,q}_p(AB)$. Identify $Q$ with
$(i,m,n)$ and let $f=\sum_{Q\in\mathcal{Q}_{AB}} s_Q \psi_Q$. By
Lemma \ref{l:Overlap_bnd}, Remark \ref{r:l:Overlap_bnd1} and Lemma
\ref{l:c:l:conv_shrlts},
$$\norm{\psi_{A^{-j}B^{[-\ell]}}\ast f}_{L^p}\leq
    C_p \sum_{i=j-1}^{j+1}\sum_{m(\ell,i)} \left(\int_{\mathbb{R}^d}
    (\sum_{Q\in\mathcal{Q}^{i,m}} \frac{\abs{s_Q}\abs{Q}^{-\frac{1}{2}}}{(1+2^j\abs{x-x_Q})^N})^p dx\right)^{1/p},$$
for all $N>d$. By the $p$-triangular inequality $\abs{a+b}^p\leq
\abs{a}^p+\abs{b}^p$ if $0<p\leq 1$ or by Hölder's inequality if
$1<p<\infty$, for a sufficiently large $N$ we have
$$\norm{\psi_{A^{-j}B^{[-\ell]}}\ast f}_{L^p}\leq
C_p \sum_{i=j-1}^{j+1}\sum_{m(\ell,i)} \left(\int_{\mathbb{R}^d}
    \sum_{Q\in\mathcal{Q}^{i,m}} \frac{\abs{s_Q}^p\abs{Q}^{-\frac{p}{2}}}{(1+2^j\abs{x-x_Q})^{d+1}})^p dx\right)^{1/p}.$$
Therefore, (notice that, since $\abs{i-j}\leq 1$, the change from
$Q\in\mathcal{Q}^{i,m}$ to $Q\in\mathcal{Q}^{j,\ell}$ is harmless)
\begin{eqnarray*}
  \norm{f}_{\mathbf{B}^{\alpha,q}_p(AB)}
    &\leq& C_{p,q} \left(\sum_{\mathfrak{d}=1}^d\sum_{j\geq 0}\sum_{\abs{[\ell]}\preceq 2^j}
        \abs{Q_j}^{-\alpha q}\left[\sum_{Q\in\mathcal{Q}^{j,\ell}} \abs{s_Q}^p\abs{Q}^{-\frac{p}{2}+\frac{d}{(d+1)}}\right]^{\frac{q}{p}}\right)^{1/q} \\
    &=& C_{p,q} \norm{\mathbf{s}}_{\mathbf{b}^{\alpha,q}_p(AB)},
\end{eqnarray*}
which is what we wanted to prove.
\hfill $\blacksquare$ \vskip .5cm   

\begin{rem}\label{r:Def_BsvSpcs_Indep_psi}
With the same arguments as in Remark 2.6 in \cite{FJ90}, the
definition of $\mathbf{B}^{\alpha,q}_p(AB)$ is independent of the
choice of $\psi\in\mathcal{S}$ as long as it satisfies the
requirements in Section \ref{S:Shrlts_HD}.
\end{rem}

\vskip 1cm
\section{Shear anisotropic inhomogeneous Triebel-Lizorkin spaces}\label{S:AB-T-L}
After defining the spaces we will leave aside the analysis of the
low frequencies and the analysis for all directions. We follow
\cite{FJ90}.

\begin{defi}\label{d:T-L_func_spcs}
Let $\Psi, \psi\in \mathcal{S}$ be as in Section \ref{S:Shrlts_HD}.
Let $\alpha\in\mathbb{R}$, $0<p<\infty$ and $0<q\leq \infty$. The
\textbf{shear anisotropic inhomogeneous Triebel-Lizorkin
\emph{distribution} space} $\mathbf{F}^{\alpha,q}_p(AB)$ is defined
as the collection of all $f\in\mathcal{S}'$ such that
\begin{eqnarray}\label{e:def_T-L_func_spcs}
\nonumber
  \norm{f}_{\mathbf{F}^{\alpha,q}_p(AB)}
    &=& \norm{f\ast \Psi}_{L^p} \\
    &+& \norm{\left(\sum_{\mathfrak{d}=1}^d\sum_{j\geq 0}\sum_{\ell=-2^j}^{2^j}
        [\abs{Q_j}^{-\alpha}\abs{\tilde{\psi}^{\mathfrak{d}}_{A^{-j}_{(\mathfrak{d})}B^{-\ell}_{(\mathfrak{d})}}\ast
        f}]^q\right)^{1/q}}_{L^p}<\infty.
\end{eqnarray}
\end{defi}

\begin{defi}\label{d:T-L_seq_spcs}
Let $\alpha\in\mathbb{R}$, $0<p<\infty$ and $0<q\leq \infty$.  The
\textbf{shear anisotropic inhomogeneous Triebel-Lizorkin
\emph{sequence} space} $\mathbf{f}^{\alpha,q}_p(AB)$ is defined as
the collection of all complex-valued sequences
$\mathbf{s}=\{s_Q\}_{Q\in\mathcal{Q}_{AB}}$ such that
\begin{equation}\label{e:def_T-L_seq_spcs}
  \norm{\mathbf{s}}_{\mathbf{f}^{\alpha,q}_p(AB)}
    = \norm{\left(\sum_{Q\in\mathcal{Q}_{AB}}(\abs{Q}^{-\alpha}\abs{s_Q}\tilde{\chi}_Q)^q\right)^{1/q}}_{L^p}<\infty.
\end{equation}
\end{defi}

\vskip0.5cm
\subsection{Two basic results}\label{sS:Two_bsc_reslts}
The proof of the characterization follows the corresponding result
in \cite{FJ90}. This is based on a kind of Peetre's inequality to
bound $S_\psi:\mathbf{F}^{\alpha,q}_p(AB)\rightarrow
\mathbf{f}^{\alpha,q}_p(AB)$, and a characterization of
$\textbf{f}_p^{\alpha,q}(AB)$ to bound $T_\psi:
\mathbf{f}^{\alpha,q}_p(AB)\rightarrow \mathbf{F}^{\alpha,q}_p(AB)$.

For all $\lambda>0$, let
\begin{equation}\label{e:def_sup_of_conv}
(\psi^{\ast\ast}_{j,\ell,\lambda}f)(x):=\sup_{y\in\mathbb{R}^d}
\frac{\abs{(\psi_{A^{-j}B^{-\ell}}\ast f)(x-y)}}{(1+\abs{B^\ell
A^jy})^{d\lambda}},
\end{equation}
be the \textbf{\emph{shear anisotropic} Peetre's maximal function}.
The next result can be regarded as a \textbf{\emph{shear
anisotropic} Peetre's inequality} and is proved in Section
\ref{sS:Proofs_AB-T-L}.

\begin{lem}\label{l:SupConv_leq_HLMaxFcn} Let $\psi$ be band
limited, $f\in \mathcal{S}'$, $j\geq 0$ and $\abs{[\ell]}\preceq
2^j$. Then, for any real $\lambda>0$, there exists a constant
$C_\lambda$ such that
$$(\psi^{\ast\ast}_{j,\ell,\lambda}f)(x)
\leq C_\lambda\left\{\mathcal{M}(\abs{\psi_{A^{-j}B^{-\ell}}\ast
f}^{1/\lambda})(x)\right\}^\lambda, \;\;\; x\in\mathbb{R}^d.$$
\end{lem}

Identify $Q$ and $P$ with $(i,m,n)$ and $(j,\ell,k)$, respectively.
For all $r>0$, $N\in\mathbb{N}$ and $i\geq j\geq 0$, define
$$(s_{r,N}^\ast)_Q:=
\left(\sum_{P\in\mathcal{Q}^{j,\ell}}\frac{\abs{s_P}^r}{(1+2^j\abs{x_Q-x_P})^N}\right)^{1/r},$$
and
$\mathbf{s}^\ast_{r,N}=\{(s^\ast_{r,N})_Q\}_{Q\in\mathcal{Q}_{AB}}$.
We then have the characterization of the sequence spaces
$\mathbf{f}^{\alpha,q}_p(AB)$ in terms of $\mathbf{s}^\ast_{r,N}$.
Next result is also proved in Section \ref{sS:Proofs_AB-T-L}.
\begin{lem}\label{l:s_ast_bnd_s}
Let $\alpha\in \mathbb{R}$, $0<p<\infty$ and $0<q\leq\infty$. Then,
for all $r>0$ and $N>(d+1)\max(1,r/q,r/p)$ there exists $C_{r,d}>0$
such that
$$\norm{\mathbf{s}}_{\mathbf{f}^{\alpha,q}_p(AB)}\leq\norm{\mathbf{s}^\ast_{r,N}}_{\mathbf{f}^{\alpha,q}_p(AB)}
\leq C_{r,d} \norm{\mathbf{s}}_{\mathbf{f}^{\alpha,q}_p(AB)}.$$
\end{lem}

\vskip0.5cm
\subsection{The characterization}\label{sS:S-T_psi-ops_Bnd}
We prove the boundedness of the analysis and synthesis operators
(\ref{e:AnlsOpShrlt}) and (\ref{e:SnhtsOpShrlt}) on the spaces in
Definitions \ref{d:T-L_func_spcs} and \ref{d:T-L_seq_spcs}.

\begin{thm}\label{t:S-T-psi-ops_Bnd} Let $\alpha\in\mathbb{R}$,
$0<p<\infty$ and $0<q\leq\infty$. Then, the operators $S_\psi :
\mathbf{F}^{\alpha,q}_p(AB)\rightarrow \mathbf{f}^{\alpha,q}_p(AB)$
and $T_\psi : \mathbf{f}^{\alpha,q}_p(AB)\rightarrow
\mathbf{F}^{\alpha,q}_p(AB)$ are well defined and bounded.
\end{thm}
\textbf{Proof.} We prove only the case $q<\infty$. To prove the
boundedness of $S_\psi$ suppose $f\in \mathbf{F}^{\alpha,q}_p(AB)$.
Let $P$ be identified with $(j,\ell,k)$. Then,
$\abs{\tilde{\psi}_{A^{-j}B^{-\ell}}\ast
f(x_P)}\chi_P=\abs{\ip{f}{\psi_P}}\tilde{\chi}_P$, as in
(\ref{e:InnrProd_equiv_Conv_Shrlts}). Let $E$ be the set of
parallelograms in $\mathcal{Q}^{j,\ell}$ surrounding the origin.
Since $\mathcal{Q}^{j,\ell}$ is a partition of $\mathbb{R}^d$ we
have for $x\in P$ with $P\in\mathcal{Q}^{j,\ell}$,
\begin{eqnarray*}
\!\!\!\!\!\!\!\!\!\!\!\!\!\!\!\!\!\!\!\!\!\!\!\!\!\!\!\!\!\!\!\!\!\!\!\!
    & & \!\!\!\!\!\!\!\!\!\!\!\!\!\!\!\!\!\!\!\!\!\!\!\!\!\!\!\!\!\!\!\!\!\!\!\!
        \sum_{P\in\mathcal{Q}^{j,\ell}} [\abs{P}^{-\alpha}\abs{(S_\psi
  f)_P}\tilde{\chi}_P(x)]^q \\
    &=& \abs{\text{det }A}^{j\alpha q} \sum_{P\in\mathcal{Q}^{j,\ell}}
        \left[\abs{\tilde{\psi}_{A^{-j}B^{-\ell}}\ast f(x_P)}\chi_P(x)\right]^q\\
    &\leq& \abs{ \text{det }A}^{j\alpha q} \sum_{P\in\mathcal{Q}^{j,\ell}}
        \sup_{y\in P}\abs{\tilde{\psi}_{A^{-j}B^{-\ell}}\ast
        f(y)}^q \chi_P(x)\\
    &\leq& \abs{ \text{det }A}^{j\alpha q} \sup_{z\in E}
        \abs{\tilde{\psi}_{A^{-j}B^{-\ell}}\ast f(x-z)}^q \\
    &=& \abs{ \text{det }A}^{j\alpha q} \sup_{z\in E}
        \left[\frac{\abs{\tilde{\psi}_{A^{-j}B^{-\ell}}\ast
        f(x-z)}}{(1+\abs{B^\ell A^j z})^{d/\lambda}}\right]^q (1+\abs{B^\ell A^j z})^{qd/\lambda} \\
    &\leq& \abs{\text{det }A}^{j\alpha q}
        \left[\sup_{z\in \mathbb{R}^d} \frac{\abs{\tilde{\psi}_{A^{-j}B^{-\ell}}\ast f(x-z)}}{(1+\abs{B^\ell A^j z})^{d/\lambda}}\right]^q
        (1+\text{Diam}(Q_{0,0,1}))^{qd/\lambda} \\
    &=& C_{d,q,\lambda}\abs{\text{det }A}^{j\alpha q}(\tilde{\psi}^{\ast\ast}_{j,\ell,1/\lambda} f)^q(x) \\
    &\leq& C_{d,q,\lambda}\abs{\text{det }A}^{j\alpha q} \left\{\mathcal{M}\left(\abs{\tilde{\psi}_{A^{-j}B^{-\ell}}\ast
        f}^\lambda\right)(x)\right\}^{q/\lambda},
\end{eqnarray*}
because of Lemma \ref{l:SupConv_leq_HLMaxFcn} (with $1/\lambda$
instead of $\lambda$ in the last inequality). Now, take
$0<\lambda<\min(p,q)$. Then, the previous estimate and Theorem
\ref{t:Feff-Stn_Ineq} yield
\begin{eqnarray*}
  \norm{S_\psi f}_{\mathbf{f}^{\alpha,q}_p(AB)}
    &=& \norm{\left(\sum_{j\geq 0}\sum_{[\ell]\preceq2^j}\sum_{P\in \mathcal{Q}^{j,\ell}}
        [\abs{P}^{-\alpha}\abs{(S_\psi f)_P}\tilde{\chi}_P]^q\right)^{1/q}}_{L^p} \\
    &\leq& C_{d,q,\lambda}\norm{\left(\sum_{j\geq 0}\sum_{{[\ell]\preceq2^j}}
        \left\{\mathcal{M}\left(\abs{\text{det }A}^{j\alpha\lambda}\abs{\tilde{\psi}_{A^{-j}B^{-\ell}}\ast
        f}^\lambda
        \right)\right\}^{q/\lambda}\right)^{1/q}}_{L^p} \\
    &=& C_{d,q,\lambda}\norm{\left(\sum_{j\geq 0}\sum_{{[\ell]\preceq2^j}}
        \left\{\mathcal{M}\left(\abs{\text{det }A}^{j\alpha\lambda}\abs{\tilde{\psi}_{A^{-j}B^{-\ell}}\ast
        f}^\lambda
        \right)\right\}^{q/\lambda}\right)^{\lambda/q}}_{L^{p/\lambda}}^{1/\lambda} \\
    &\leq& C_{d,p,q,\lambda}\norm{\left(\sum_{j\geq 0}\sum_{{[\ell]\preceq2^j}}
        \abs{\text{det }A}^{j\alpha q}\abs{\tilde{\psi}_{A^{-j}B^{-\ell}}\ast f}^q
        \right)^{\lambda/q}}_{L^{p/\lambda}}^{1/\lambda} \\
    &=& C_{d,p,q,\lambda}\norm{\left(\sum_{j\geq 0}\sum_{{[\ell]\preceq2^j}}
        [\abs{\text{det }A}^{j\alpha}\abs{\tilde{\psi}_{A^{-j}B^{-\ell}}\ast
        f}]^q
        \right)^{1/q}}_{L^{p}} \\
    &=& C_{d,p,q,\lambda}\norm{f}_{\mathbf{F}^{\alpha,q}_p(AB)}.
\end{eqnarray*}

To prove the boundedness of $T_\psi$ suppose
$\mathbf{s}=\{s_Q\}_{Q\in\mathcal{Q}_{AB}}\in
\mathbf{f}^{\alpha,q}_p$ and $f=T_\psi\mathbf{s}=\sum_{Q\in
\mathcal{Q}_{AB}}s_Q\psi_Q$. Identify $Q$ with $(i,m,n)$. By Lemma
\ref{l:Overlap_bnd} (see also Remark \ref{r:l:Overlap_bnd1}) and
Lemma \ref{l:c:l:conv_shrlts}, we have for $x\in Q$ with
$Q\in\mathcal{Q}^{i,m}$,
\begin{eqnarray*}
  \abs{\tilde{\psi}_{A^{-j}B^{-\ell}}\ast f(x)}
    &\leq& \sum_{i=j-1}^{j+1}\sum_{m(\ell,i)}\sum_{Q\in\mathcal{Q}^{i,m}}
        \abs{s_Q}\abs{\tilde{\psi}_{A^{-j}B^{-\ell}}\ast \psi_Q(x)} \\
    &\leq& C_N\sum_{i=j-1}^{j+1}\sum_{m(\ell,i)}\sum_{Q\in\mathcal{Q}^{i,m}}
        \abs{s_Q}\frac{\abs{Q}^{-1/2}}{(1+2^i\abs{x-x_Q})^N} \\
    &\leq& C_N'\sum_{i=j-1}^{j+1}\sum_{m(\ell,i)}\sum_{Q\in\mathcal{Q}^{i,m}}
        \abs{s_Q}\frac{\abs{Q}^{-1/2}}{(1+2^i\abs{x_{Q'}-x_Q})^N} \\
    &=& C_N'\sum_{i=j-1}^{j+1}\sum_{m(\ell,i)}
        \abs{Q}^{-1/2}(s_{1,N}^\ast)_{Q}\chi_{Q}(x)\\
    &=& C_N'\sum_{i=j-1}^{j+1}\sum_{m(\ell,i)}\sum_{Q\in\mathcal{Q}^{i,m}}
        (s_{1,N}^\ast)_Q \tilde{\chi}_Q(x),
\end{eqnarray*}
for all $N>d$ and because $\mathcal{Q}^{i,m}$ is a partition of
$\mathbb{R}^d$. Let $N>(d+1)\max(1,1/q,1/p)$. Then, (notice that,
since $\abs{i-j}\leq 1$, the change $Q\in\mathcal{Q}^{i,m}$ for
$Q\in\mathcal{Q}^{j,\ell}$ is harmless) the previous estimate yields
\begin{eqnarray*}
        \norm{T_\psi \mathbf{s}}_{\mathbf{F}^{\alpha,q}_p(AB)}
    &\leq& C_{d,p,q}\norm{\left(\sum_{j\geq 0}\sum_{{[\ell]\preceq2^j}}
        \left[\sum_{Q\in\mathcal{Q}^{j,\ell}}
        \abs{Q}^{-\alpha}(s_{1,N}^\ast)_Q\tilde{\chi}_Q\right]^q\right)^{1/q}}_{L^p} \\
    &=& C_{d,p,q} \norm{\left(\sum_{j\geq 0}\sum_{{[\ell]\preceq2^j}}
        \sum_{Q\in\mathcal{Q}^{j,\ell}}[
        \abs{Q}^{-\alpha}(s_{1,N}^\ast)_Q\tilde{\chi}_Q]^q\right)^{1/q}}_{L^p}\\
    &=& C_{d,p,q}\norm{\mathbf{s}_{1,N}^\ast}_{\mathbf{f}^{\alpha,q}_p(AB)}
        \leq
        C_{d,p,q}\norm{\mathbf{s}}_{\mathbf{f}^{\alpha,q}_p(AB)},
\end{eqnarray*}
because $\mathcal{Q}^{j,\ell}$ is a partition of $\mathbb{R}^d$ and
Lemma \ref{l:s_ast_bnd_s} in the last inequality.
\hfill $\blacksquare$ \vskip .5cm   

\begin{rem}\label{r:Def_T-LSpcs_Indep_psi}
With the same arguments as in Remark 2.6 in \cite{FJ90}, the
definition of $\mathbf{F}^{\alpha,q}_p(AB)$ is independent of the
choice of $\psi\in\mathcal{S}$ as long as it satisfies the
requirements in Section \ref{S:Shrlts_HD}.
\end{rem}

\vskip 1cm \section{A reproducing identity with smooth Parseval
frames}\label{S:ReprId_SmthShrlts} Recently, Guo and Labate in
\cite{GL11} found a way to overcome the use of characteristic
functions in the Fourier domain to restrict the shearlets to the
respective cone (see (\ref{e:ShrltPrsvlFrm})). The use of these
characteristic functions affects the smoothness of the
\emph{boundary shearlets} (those with $\abs{[\ell]}=\pm 2^j$) in the
Fourier domain and, therefore, their spatial localization. They
slightly modify the definition of these boundary shearlets instead
of projecting them into the cone. This new shearlets system is not
affine-like. However, they do produce the same frequency covering as
that in Section \ref{S:Shrlts_HD}.

\vskip 0.5cm
\subsection{The new smooth shearlets system}\label{sS:SmthPrsvlFrms}
This subsection is a brief summary of some results in \cite{GL11}
and is intended to show the construction of such smooth Parseval
frames. Let $\hat{\phi}$ be a $C^\infty$ univariate function such
that $0\leq\hat{\phi}\leq 1$, with $\hat{\phi}=1$ on $[-1/16,1/16]$
and $\hat{\phi}=0$ outside $[-1/8,1/8]$ (\emph{i.e.}, $\phi$ is a
rescaled scaling function of a Meyer wavelet). For
$\xi\in\hat{\mathbb{R}}^d$, let
$\hat{\Phi}(\xi)=\hat{\phi}(\xi_1)\hat{\phi}(\xi_2)\cdots\hat{\phi}(\xi_d)$
and $W^2(\xi)=\hat{\Phi}^2(2^{-2}\xi)-\hat{\Phi}^2(\xi)$. It follows
that
$$\hat{\Phi}(\xi)+\sum_{j\geq 0}W^2(2^{-2j}\xi)=1, \;\;\text{ for all } \xi\in\hat{\mathbb{R}}^d.$$
Let now $v\in C^\infty(\mathbb{R})$ be such that $v(0)=1$,
$v^{(n)}(0)=0$ for all $n\geq 1$, supp $v\subset [-1,1]$ and
$$\abs{v(u-1)}^2+\abs{v(u)}^2+\abs{v(u+1)}^2=1, \;\; \abs{u}\leq 1.$$
Then, for any $j\geq 0$,
$$\sum_{m=-2^j}^{2^j} \abs{v(2^ju-m)}^2=1, \;\;\abs{u}\leq 1.$$
See \cite{GL11} for comments on the construction of these functions
and their properties. With
$V_\frak{d}(\xi)=\prod_{i\not=\frak{d}}v(\frac{\xi_i}{\xi_\frak{d}})$,
$\xi\in\mathcal{D}^{(\frak{d})}$, the shearlets system for
$L^2((\mathcal{D}^{(\frak{d})})^\vee)$ is defined as the countable
collection of functions
$$\{\psi^{(\frak{d})}_{j,\ell,k}: \frak{d}=1,\ldots,d, j\geq 0, \abs{[\ell]}\preceq2^j, k\in\mathbb{Z}^d\},$$
whose ``inner" elements ($\abs{[\ell]}\prec2^j$) are defined by
their Fourier transform
\begin{equation}\label{e:def_new_Shrlt_Sys_Fourier}
(\psi^{(\frak{d})}_{j,\ell,k})^\wedge(\xi)=\abs{\text{det }
A_{(\frak{d})}}^{-j/2}
    W(2^{-2j}\xi)V_\frak{d}(\xi A_{(\frak{d})}^{-j}B_{(\frak{d})}^{[-\ell]})\mathbf{e}^{-2\pi i \xi
    A_{(\frak{d})}^{-j}B_{(\frak{d})}^{[-\ell]}k}, \;\;\; \xi\in\mathcal{D}^{(\frak{d})}.
\end{equation}
The \emph{boundary shearlets} are defined slightly different but
share similar properties in both the time and Fourier domain. This
new system is not affine-like since the function $W$ is not
shear-invariant. However, they generate the same covering of
$\hat{\mathbb{R}}^d$. The new smooth Parseval frame condition is now
written as (see Theorem 2.3 in \cite{GL11})
\begin{equation}\label{e:SmthPrsvlFrm_Prop_Shrlts}
\abs{\hat{\Psi}(\xi)}^2
    +\sum_{\mathfrak{d}=1}^d\sum_{j\geq0}\sum_{\abs{[\ell]}\prec2^j}
        \abs{\hat{\psi}^{(\mathfrak{d})}(\xi A_{(\mathfrak{d})}^{-j} B_{(\mathfrak{d})}^{[-\ell]})}^2
    + \sum_{\mathfrak{d}=1}^d\sum_{j\geq 0}\sum_{\abs{[\ell]}=\pm 2^j}
        \abs{\hat{\psi}^{(\mathfrak{d})}(\xi A_{(\mathfrak{d})}^{-j} B_{(\mathfrak{d})}^{[-\ell]})}^2=1,
\end{equation}
for all $\xi\in\hat{\mathbb{R}}^d$. Notice that now there do not
exist characteristic functions as in (\ref{e:ShrltPrsvlFrm}).

\vskip0.5cm
\subsection{The reproducing identity on $\mathcal{S}'$}\label{sS:Repr_Ident_Distrbtns}
Our goal is to show that, with the smooth Parseval frames of
shearlets of Guo and Labate in \cite{GL11}, $T_\psi\circ S_\psi$ is
the identity on $\mathcal{S}'$. First, we show that any
$f\in\mathcal{S}'$ admits a kind of Littlewood-Paley decomposition
with shear anisotropic dilations, for which we follow \cite{BH05}.
Then, we show the reproducing identity in $\mathcal{S}'$ following
\cite{FJW}.

\begin{lem}\label{l:Ident_Conv}
Let $\{\Psi(\cdot-k)\}_{k\in\mathbb{Z}^d}\cup\{\psi_{j,\ell,k}:j\geq
0, \abs{[\ell]}\preceq 2^j, k\in\mathbb{Z}^d\}$ be the smooth
shearlet system that verifies (\ref{e:SmthPrsvlFrm_Prop_Shrlts}).
Then, for any $f\in\mathcal{S}'$,
\begin{eqnarray*}
  f=f\ast\tilde{\Psi}\ast\Psi
    &+& \sum_{\mathfrak{d}=1}^d\sum_{j\geq 0}\sum_{\abs{[\ell]}\prec2^j}
        f\ast\tilde{\psi}^{(\mathfrak{d})}_{A_{(\mathfrak{d})}^{-j}B_{(\mathfrak{d})}^{[-\ell]}}
            \ast\psi^{(\mathfrak{d})}_{A_{(\mathfrak{d})}^{-j}B_{(\mathfrak{d})}^{[-\ell]}} \\
    &+& \sum_{\mathfrak{d}=1}^d\sum_{j\geq 0}\sum_{\abs{[\ell]}= 2^j}
        f\ast\tilde{\psi}^{(\frak{d})}_{A_{(\mathfrak{d})}^{-j}B_{(\mathfrak{d})}^{[-\ell]}}
            \ast\psi^{(\frak{d})}_{A_{(\mathfrak{d})}^{-j}B_{(\mathfrak{d})}^{[-\ell]}},
\end{eqnarray*}
with convergence in $\mathcal{S}'$.
\end{lem}
\textbf{Proof.} One can see Peetre's discussion on pp. 52-54 in
\cite{Pe76} regarding convergence. Since the Fourier transform
$\mathcal{F}$ is an isomorphism of $\mathcal{S}'$, it suffices to
show that
\begin{eqnarray*}
  \hat{f}(\xi)=\hat{f}(\xi) \abs{\hat{\Psi}(\xi)}^2
    &+& \sum_{\mathfrak{d}=1}^d\sum_{j\geq 0}\sum_{\abs{[\ell]}\prec2^j}
        \hat{f}(\xi) \abs{\hat{\psi}^{(\mathfrak{d})}(\xi A_{(\mathfrak{d})}^{-j} B_{(\mathfrak{d})}^{-\ell})}^2\\
    &+& \sum_{\mathfrak{d}=1}^d\sum_{j\geq 0}\sum_{\abs{[\ell]}=\pm 2^j}
        \hat{f}(\xi) \abs{\hat{\psi}^{(\mathfrak{d})}(\xi A_{(\mathfrak{d})}^{-j} B_{(\mathfrak{d})}^{-\ell})}^2
\end{eqnarray*}
converges in $\mathcal{S}'$. Since the equality is a straight
consequence of (\ref{e:SmthPrsvlFrm_Prop_Shrlts}), we will only show
convergence in $\mathcal{S}'$ of the right-hand side of the equality
for those shearlets with $j\geq 0$ ($\Psi$ is in fact a scaling
function of a Meyer wavelet). Suppose that $\hat{f}$ has order $\leq
m$. This is, there exists an integer $n\geq 0$ and a constant $C$
such that
$$\abs{\ip{\hat{f}}{g}}\leq C \sup_{\abs{\alpha}\leq n, \abs{\beta}\leq m} \norm{g}_{\alpha,\beta}, \;\;\; \text{ for all } g\in \mathcal{S},$$
where $\norm{g}_{\alpha,\beta}=
\sup_{\xi\in\hat{\mathbb{R}}^d}\abs{\xi^\alpha}\abs{\partial^\beta
g(\xi)}$ denotes the usual semi-norm in $\mathcal{S}$ for
multi-indices $\alpha$ and $\beta$. Then,
$$\abs{\ip{\hat{f}\abs{(\psi_{A^{-j}B^{-\ell}})^\wedge}^2}{g}}=\abs{\ip{\hat{f}}{\abs{(\psi_{A^{-j}B^{-\ell}})^\wedge}^2 g}}
    \leq C \sup_{\abs{\alpha}\leq n, \abs{\beta}\leq m} \norm{\abs{(\psi_{A^{-j}B^{-\ell}})^\wedge}^2 g}_{\alpha,\beta}.$$
As in Lemma 2.5 in \cite{GL07}, one can prove that
$$\sup_{\abs{\beta}=m} \norm{\partial^\beta \abs{(\psi_{A^{-j}B^{-\ell}})^\wedge}^2}_\infty \leq C 2^{-jm}.$$
Hence, by the compact support conditions of
$(\psi_{A^{-j}B^{-\ell}})^\wedge(\xi)$ (see Section \ref{S:Intro})
\begin{eqnarray*}
    & & \!\!\!\!\!\!\!\!\!\!\!\!\!\!\!\!\!\!\!\!\!\!\!\!\!\!\!\!\!\!\!\!\!\!\!\!\!\!\!\!
            \sup_{\abs{\alpha}\leq n, \abs{\beta}\leq m} \norm{\abs{(\psi_{A^{-j}B^{-\ell}})^\wedge}^2
            g}_{\alpha,\beta}\\
    &\leq& C \sup_{\xi\in \hat{\mathbb{R}}^2} \left[(1+\abs{\xi})^n
            \sup_{\abs{\beta}\leq m} \abs{\partial^\beta\abs{(\psi_{A^{-j}B^{-\ell}})^\wedge(\xi)}^2}
            \sup_{\abs{\beta}\leq m} \abs{\partial^\beta g(\xi)}\right] \\
    &\leq& C \sup_{\xi\in \text{supp}(\psi_{A^{-j}B^{-\ell}})^\wedge(\xi)} (1+\abs{\xi})^n \sup_{\abs{\beta}\leq m}\abs{\partial^\beta g(\xi)} \\
    &\leq& C \sup_{\abs{\alpha}\leq n+1, \abs{\beta}\leq m} \norm{g}_{\alpha,\beta}
            \sup_{\xi\in\text{supp}(\psi_{A^{-j}B^{-\ell}})^\wedge(\xi)} (1+\abs{\xi})^{-1}\\
    &\leq& C \sup_{\abs{\alpha}\leq n+1, \abs{\beta}\leq m} \norm{g}_{\alpha,\beta}
            (1+2^{2j-4})^{-1}\leq C2^{-2j},
\end{eqnarray*}
which proves the convergence in $\mathcal{S}'$.
\hfill $\blacksquare$ \vskip .5cm   

\begin{thm}\label{t:Ident_in_Distributions}
Let the shearlet system $\{\psi_{j,\ell,k}\}$ be constructed as in
Subsection \ref{sS:SmthPrsvlFrms} such that it is a smooth Parseval
frame that verifies (\ref{e:SmthPrsvlFrm_Prop_Shrlts}). The
composition of the analysis and synthesis operators $T_\psi\circ
S_\psi$ (see (\ref{e:AnlsOpShrlt}) and (\ref{e:SnhtsOpShrlt}) for
the definitions) is the identity
$$f=\sum_{Q\in\mathcal{Q}_{AB}} \ip{f}{\psi_Q}\psi_Q,$$
in $\mathcal{S}'$.

\end{thm}
\textbf{Proof.} As in (\ref{e:InnrProd_equiv_Conv_Shrlts}),
$f\ast\tilde{\psi}_{A^{-j}B^{-\ell}}(A^{-j}B^{-\ell}k)=
f\ast\tilde{\psi}_{A^{-j}B^{-\ell}}(x_P) = \abs{\text{det
}A}^{j/2}\ip{f}{\psi_P}= \abs{P_j^{-1/2}}\ip{f}{\psi_P}$, where $P$
is identified with $(j,\ell,k)$. Let
$g=f\ast\tilde{\psi}_{A^{-j}B^{-\ell}}$ and
$h=\psi_{A^{-j}B^{-\ell}}$. By construction, 
supp$(\psi_{j,\ell,k})^\wedge(\xi)\subset QB^\ell A^j$. Therefore,
Lemma \ref{l:Smplng_BndLim_Shrlt} yields
\begin{eqnarray*}
  f\ast\tilde{\psi}_{A^{-j}B^{-\ell}}\ast\psi_{A^{-j}B^{-\ell}}
    &=& \sum_{k\in\mathbb{Z}^d} \ip{f}{\psi_{j,\ell,k}}\psi_{j,\ell,k} \\
    &=& \sum_{P\in\mathcal{Q}^{j,\ell}} \ip{f}{\psi_P}\psi_P.
\end{eqnarray*}
By appropriately summing over $\mathfrak{d}=1,\ldots,d$, $j\geq 0$
and $\abs{[\ell]}\preceq 2^j$, Lemma \ref{l:Ident_Conv} yields the
result.
\hfill $\blacksquare$ \vskip .5cm   

\vskip 1cm
\section{Embeddings}\label{S:Embeddings}

\vskip 0.5cm
\subsection{Sobolev-type embeddings}\label{sS:Sblv-type_Emmbeddings}
To prove the next embeddings we follow Section 2.3.2 in \cite{Tr83}.
\begin{thm}\label{t:Sobolev-type_embeddings}
Let $s\in\mathbb{R}$.
\begin{itemize}
  \item[i)] For $0<q_1\leq q_0\leq \infty$,
        $$\mathbf{B}^{s,q_1}_p(AB)\hookrightarrow \mathbf{B}^{s,q_0}_p(AB), \;\;\; 0<p\leq\infty$$
        and
        $$\mathbf{F}^{s,q_1}_p(AB)\hookrightarrow \mathbf{F}^{s,q_0}_p(AB), \;\;\; 0<p<\infty.$$
  \item[ii)] For $0< q_0\leq \infty$, $0< q_1\leq \infty$ and
        $\varepsilon>\frac{d-1}{(d+1)q_0}$,
        $$\mathbf{B}^{s+\varepsilon,q_1}_p(AB)\hookrightarrow \mathbf{B}^{s,q_0}_p(AB), \;\;\; 0<p\leq\infty$$
        and
        $$\mathbf{F}^{s+\varepsilon,q_1}_p(AB)\hookrightarrow \mathbf{F}^{s,q_0}_p(AB), \;\;\; 0<p<\infty.$$
  \item[iii)] For $0<q\leq\infty$, $0<p<\infty$ and
        $s\in\mathbb{R}$,
        $$\mathbf{B}^{s,\min{\{p,q\}}}_p(AB)\hookrightarrow \mathbf{F}^{s,q}_p(AB) \hookrightarrow \mathbf{B}^{s,\max{\{p,q\}}}_p(AB).$$
\end{itemize}
\end{thm}
\textbf{Proof}. The monotonicity of $\ell^q$ norms proves i). To
prove ii) let $\varepsilon>\frac{d-1}{(d+1)q_0}$. Then, the Besov
case follows from
\begin{eqnarray*}
  \norm{f}_{\mathbf{B}^{s,q_0}_p}
    &\leq& \sup_{\frak{d},j,[\ell]} \abs{Q_j}^{-(s+\varepsilon)}\norm{f\ast \psi^{(\frak{d})}_{A^{-j}B^{-\ell}}}_{L^p}
        \left(\sum_{\frak{d}',j',[\ell ']} \abs{Q_j}^{\varepsilon q_0}\right)^{\frac{1}{q_0}} \\
    &\lesssim& \sup_{\frak{d},j,[\ell]} \abs{Q_j}^{-(s+\varepsilon)}\norm{f\ast \psi^{(\frak{d})}_{A^{-j}B^{-\ell}}}_{L^p}
        \left(\sum_{j'\geq 0} 2^{-j((d+1) \varepsilon q_0-(d-1))}\right)^{\frac{1}{q_0}} \\
    &=& C_{\varepsilon,q_0,d} \sup_{\frak{d},j,[\ell]} \abs{Q_j}^{-(s+\varepsilon)}\norm{f\ast
    \psi^{(\frak{d})}_{A^{-j}B^{-\ell}}}_{L^p},
\end{eqnarray*}
and $\ell^\infty\hookrightarrow\ell^{q_1}$. Similarly for the
Triebel-Lizorkin case in ii). To prove iii) write
$a_{j,\ell}(x)=\abs{Q_j}^{-\alpha}f\ast \psi_{A^{-j} B^{-\ell}}(x)$.
Since the set $\{j\geq 0, [\abs{\ell}]\preceq 2^j\}$ is countable,
one can find a bijection with $k\in\mathbb{Z}$. Consider first
$0<q\leq p<\infty$. Then,
\begin{eqnarray*}
  \left(\sum_k \norm{a_k}^p_{L^p}\right)^{\frac{1}{p}}
    &\leq& \left(\int_{\mathbb{R}^n} (\sum_k \abs{a_k(x)}^q)^{\frac{p}{q}}dx\right)^{\frac{1}{p}} \\
    &\leq& \left(\sum_k \norm{\abs{a_k(x)}^q}_{L^{p/q}}\right)^{\frac{1}{q}}
        = \left(\sum_k \norm{a_k}_{L^p}^q\right)^{\frac{1}{q}},
\end{eqnarray*}
because of usual $\ell^p$ inequalities and Minkowski's inequality.
Now let $0<p<q \leq\infty$. Then,
\begin{eqnarray*}
  \left(\sum_k \norm{a_k}_{L^p}^q\right)^{\frac{1}{q}}
    &=& \norm{\int \abs{a_k(x)}^p dx}_{\ell^{q/p}}^{\frac{1}{p}}
        \leq \left(\int \norm{\abs{a_k(x)}^p}_{\ell^{q/p}} dx\right)^{\frac{1}{p}}\\
    &=& \left(\int \norm{\abs{a_k(x)}}^p_{\ell^q} dx\right)^{\frac{1}{p}}
        \leq \left(\int \norm{a_k(x)}_{\ell^p}^p
        dx\right)^{\frac{1}{p}},
\end{eqnarray*}
because of the generalized Minkowski's inequality and usual $\ell^p$
inequalities.

\hfill $\blacksquare$ \vskip .5cm   

\vskip 0.5cm
\subsection{Embeddings of Besov spaces}\label{sS:Emmbeddings_Besov} We now
present embeddings between $\mathbf{B}^{\alpha_1,q}_p$ and
$\mathbf{B}^{\alpha_2,q}_p(AB)$, for certain conditions on the
smoothness parameters $\alpha_1$ and $\alpha_2$.

\begin{thm}\label{t:Bsv_in_BsvAB}
Let $0<p,q\leq \infty$ and $\alpha_1,\alpha_2\in\mathbb{R}$. For
$0<p\leq 1$ let $\lambda=-1$ and for $1<p<\infty$ let $\lambda>p/2$
¿CASO $p=\infty $?. Then,
$$\mathbf{B}^{\alpha_1,q}_p\hookrightarrow
\mathbf{B}^{\alpha_2,q}_p(AB),$$ when
$\frac{2d\lambda}{p}+\frac{d-1}{q}+(d+1)\alpha_2<2\alpha_1$.
\end{thm}
\textbf{Proof.} To shorten notation write
$\mathbf{b}_1=\mathbf{b}^{\alpha_1,q}_p$,
$\mathbf{B}_1=\mathbf{B}^{\alpha_1,q}_p$ and
$\mathbf{B}_2=\mathbf{B}^{\alpha_2,q}_p(AB)$. Let
$f=\sum_{Q\in\mathcal{D}_+}s_Q\varphi_Q \in\mathbf{B}_1$. From the
compact support conditions on $(\varphi_{\nu,k})^\wedge$ and
$(\psi_{A^{-j}B^{-\ell}})^\wedge$, we have
\begin{eqnarray*}
  \abs{f\ast\psi_{A^{-j}B^{-\ell}}(x)}^p
    &=& \abs{\sum_{\nu\geq 0}\sum_{k\in\mathbb{Z}^d} s_{\nu,k}\varphi_{\nu,k}\ast\psi_{A^{-j}B^{-\ell}}(x)}^p \\
    &\lesssim& \left(\sum_{k\in\mathbb{Z}^d}\abs{s_{2j,k}}\abs{\varphi_{2j,k}\ast\psi_{A^{-j}B^{-\ell}}(x)}\right)^p \\
    &\leq& C_N \left(\sum_{k\in\mathbb{Z}^d}\abs{s_{2j,k}}
        \frac{\abs{Q_{2j}}^{-1/2}}{(1+2^j\abs{x+2^{-2j}k})^N}\right)^p,
\end{eqnarray*}
for every $N>d$ by Lemma \ref{l:conv_shrlts-wvlts}. If $0<p\leq 1$,
choose $N>d/p$ and use the $p$-triangle inequality
$\abs{a+b}^p\leq\abs{a}^p+\abs{b}^p$ to get
\begin{eqnarray*}
  \norm{f\ast\psi_{A^{-j}B^{-\ell}}}_{L^p}
    &\leq& C_N \left(\sum_{k\in\mathbb{Z}^d}\abs{s_{2j,k}}^p\abs{Q_{2j}}^{-\frac{p}{2}}
        \int_{\mathbb{R}^d}\frac{dx}{(1+2^j\abs{x+2^{-2j}k})^{Np}}\right)^{1/p} \\
    &\leq& C_N
    \left(\sum_{k\in\mathbb{Z}^d}\abs{s_{2j,k}}^p\abs{Q_{2j}}^{p(-\frac{1}{2}+\frac{1}{2p})}\right)^{1/p}.
\end{eqnarray*}
For $1<p<\infty$ choose $N=a+b$ such that $a>d/p$ and $b>d(p-1)/p$.
Hölder's inequality yields
\begin{eqnarray*}
  \abs{f\ast\psi_{A^{-j}B^{-\ell}}(x)}^p
    &\leq& C_N \left(\sum_{k\in\mathbb{Z}^d}\frac{\abs{s_{2j,k}}\abs{Q_{2j}}^{-\frac{1}{2}}}{(1+2^j\abs{x+2^{-2j}k})^N}
        \cdot\frac{2^{jN}}{2^{jN}}\right)^p \\
    &\leq& C_N 2^{jNp}
        \left(\sum_{k\in\mathbb{Z}^d}\frac{\abs{s_{2j,k}}^p\abs{Q_{2j}}^{-\frac{p}{2}}}{(1+\abs{2^{2j}x+k})^{ap}}\right)
        \left(\sum_{k\in\mathbb{Z}^d}\frac{1}{(1+\abs{2^{2j}x+k})^{bp'}}\right)^{p/p'} \\
    &\leq& C_{N,p}
        2^{jNp}\left(\sum_{k\in\mathbb{Z}^d}\frac{\abs{s_{2j,k}}^p\abs{Q_{2j}}^{-\frac{p}{2}}}{(1+\abs{2^{2j}x+k})^{ap}}\right).
\end{eqnarray*}
With $\lambda=\frac{Np}{2d}>p/2$, we write
$2^{jNp}=\abs{Q_{2j}}^{-\frac{pN}{2d}}=\abs{Q_{2j}}^{-\lambda}$.
Since $ap>d$ we have
\begin{eqnarray*}
  \norm{f\ast\psi_{A^{-j}B^{-\ell}}}_{L^p}
    &\leq& C_{N,p} \left(\int_{\mathbb{R}^d}\sum_{k\in\mathbb{Z}^d}
        \frac{\abs{s_{2j,k}}^p\abs{Q_{2j}}^{p(-\frac{1}{2}-\frac{\lambda}{p})}}{(1+\abs{2^{2j}x+k})^{ap}} dx\right)^{1/p} \\
    &=& C_{N,p} \left(\sum_{k\in\mathbb{Z}^d}
        \abs{s_{2j,k}}^p\abs{Q_{2j}}^{p(-\frac{1}{2}-\frac{\lambda}{p}+\frac{1}{p})}\right)^{1/p}.
\end{eqnarray*}
Let $\lambda=-1$ for $0<p\leq1$ and $\lambda>p/2$ for $1<p<\infty$.
Then, since $\abs{\{\ell:\abs{[\ell]}\leq 2^j\}}\sim
2^{(d-1)j}=\abs{Q_{2j}}^{-\frac{d-1}{2d}}$ and
$\abs{P_j}^{-\alpha_2}=2^{j(d+1)\alpha_2}=\abs{Q_{2j}}^{-\frac{(d+1)\alpha_2}{2d}}$,
we have
\begin{eqnarray*}
  \norm{f}_{\mathbf{B}_2}
    &\leq& C_{N,p} \left(\sum_{j\geq 0}\sum_{\abs{[\ell]}\leq 2^j}
        \left[\abs{P_j}^{-\alpha_2}\left(\sum_{k\in\mathbb{Z}^d}
        [\abs{Q_{2j}}^{-\frac{1}{2}+\frac{1}{p}-\frac{\lambda}{p}}\abs{s_{2j,k}}]^p\right)^{1/p}\right]^q\right)^{1/q} \\
    &\leq& C_{N,p} \left(\sum_{j\geq 0}
        \left(\sum_{k\in\mathbb{Z}^d}
        [\abs{Q_{2j}}^{-\frac{1}{2}+\frac{1}{p}-\frac{\lambda}{p}-\frac{d-1}{2qd}-\frac{(d+1)\alpha_2}{2d}}\abs{s_{2j,k}}]^p\right)^{q/p}\right)^{1/q} \\
    &=& C_{N,p} \left(\sum_{j\geq 0}
        \left(\sum_{Q\in\mathcal{D}^{2j}}
        [\abs{Q}^{-\frac{1}{2}+\frac{1}{p}-\frac{\lambda}{p}-\frac{d-1}{2qd}-\frac{(d+1)\alpha_2}{2d}}\abs{s_Q}]^p\right)^{q/p}\right)^{1/q} \\
    &\leq& C_{N,p} \left(\sum_{j\geq 0}
        \left(\sum_{Q\in\mathcal{D}^{2j}}
        [\abs{Q}^{-\frac{\alpha_1}{d}+\frac{1}{p}-\frac{1}{2}}\abs{s_Q}]^p\right)^{q/p}\right)^{1/q} \\
    &\leq& C_{N,p} \left(\sum_{j\geq 0}
        \left(\sum_{Q\in\mathcal{D}^j}
        [\abs{Q}^{-\frac{\alpha_1}{d}+\frac{1}{p}-\frac{1}{2}}\abs{s_Q}]^p\right)^{q/p}\right)^{1/q}
        =C_{N,p} \norm{\mathbf{s}}_{\mathbf{b}_1}.
\end{eqnarray*}
By Theorem 2.6 in \cite{FJ85} (with slight different notation for
the analysis operator and, therefore, for the definition of the
space), the last norm is bounded by $\norm{f}_{\mathbf{B}_1}$. This
completes the proof.
\hfill $\blacksquare$ \vskip .5cm   

\begin{thm}\label{t:BsvAB_in_Bsv}
Let $0<p,q\leq \infty$, $s=[\max(1,1/p)-\min(1,1/q)]/2$ and
$\alpha_1,\alpha_2\in\mathbb{R}$. Let $\lambda=0$ if $0<p\leq 1$ or
$\lambda>d(p-1)/p$ if $1<p<\infty$ and
$2d+\lambda+2(\alpha_1+s(d-1))<(d+1)(\alpha_2+1)$. Then,
$$\mathbf{B}^{\alpha_2,q}_p(AB)\hookrightarrow \mathbf{B}^{\alpha_1,q}_p$$
\end{thm}
\textbf{Proof.} Let $f=\sum_{P\in\mathcal{Q}_{AB}}s_P\psi_P\in
\mathbf{B}^{\alpha_2,q}_p(AB)$. To shorten notation write
$\mathbf{B}_1=\mathbf{B}^{\alpha_1,q}_p$ and
$\mathbf{B}_2=\mathbf{B}^{\alpha_2,q}_p(AB)$. Let $p<1$ and $q\geq
1$. Since $\norm{\cdot}^p_{L^p}$ satisfies the triangular
inequality, using Hölder's inequality ($1/p>1$) we have that
\begin{eqnarray*}
  \norm{f}_{\mathbf{B}_1}
    &\lesssim& \left(\sum_{j=0}^{\infty} \abs{Q_{2j}}^{-\frac{\alpha_1 q}{d}}
        \left[\sum_{\abs{[\ell]}\preceq 2^j}\norm{\sum_{P\in\mathcal{Q}^{j,\ell}}
            s_P\psi_P\ast\varphi_{2^{2j} I}}^p_{L^p}\right]^{q/p}\right)^{1/q} \\
    &\lesssim& \left(\sum_{j=0}^{\infty} \abs{Q_{2j}}^{-\frac{\alpha_1 q}{d}}2^{(d-1)j q(\frac{1}{p}-1)}
        \left[\sum_{\abs{[\ell]}\preceq 2^j}\norm{\sum_{P\in\mathcal{Q}^{j,\ell}}
            s_P\psi_P\ast\varphi_{2^{2j} I}}_{L^p}\right]^q
            \right)^{1/q},
\end{eqnarray*}
since $\abs{\{\ell:\abs{[\ell]}\leq 2^j\}}\sim 2^{(d-1)j}$. Using
Hölder's inequality again ($q\geq 1$),
$$\norm{f}_{\mathbf{B}_1}\lesssim \left(\sum_{j=0}^\infty \abs{Q_{2j}}^{-\frac{\alpha_1 q}{d}}2^{2j (d-1)q(\frac{1}{p}-\frac{1}{q})/2}
    \sum_{\abs{[\ell]}\preceq 2^j}\norm{\sum_{P\in\mathcal{Q}^{j,\ell}}
    s_P\psi_P\ast\varphi_{2^{2j} I}}^q_{L^p}\right)^{1/q}.$$
The other cases are result of the triangular inequality when $p\geq
1$ and the fact that $\ell^q\hookrightarrow \ell^1$ when $q<1$.
Since $\psi_P(x)=\abs{P_j}^{-\frac{1}{2}}\psi(B^{[\ell]}A^{j}x-k)$
and $\varphi_{2^{2j}I}(x)=\abs{Q_{2j}}^{-1}\varphi(2^{2j}x)$, Lemma
\ref{l:conv_shrlts-wvlts} yields
\begin{eqnarray*}
    & &\!\!\!\!\!\!\!\!\!\!\!\!\!\!\!\!\!\!\!\!\!\!\!\!
        \norm{\sum_{P\in\mathcal{Q}^{j,\ell}}
        s_P\psi_P\ast\varphi_{2^{2j}I}}^q_{L^p}\\
    &\leq& \left(\int_{\mathbb{R}^d}\left[\sum_{P\in\mathcal{Q}^{j,\ell}}
        \abs{s_P}\abs{\psi_P\ast\varphi_{2^{2j}I}(x)}\right]^p dx\right)^{q/p} \\
    &\leq& \left(\int_{\mathbb{R}^d}\left[\sum_{k\in\mathbb{Z}^d}
        \abs{s_{j,\ell,k}}
        \frac{\abs{P_j}^{-\frac{1}{2}}\abs{Q_{2j}}^{-1}\abs{P_j}}{(1+2^j\abs{x-A^{-j}B^{-[\ell]}k})^N}
        \right]^p dx\right)^{q/p},
\end{eqnarray*}
for every $N>d$. When $0<p\leq 1$ choose $Np>d$. Then, the
$p$-triangular inequality $\abs{a+b}^p\leq\abs{a}^p+\abs{b}^p$
yields
\begin{eqnarray*}
  \norm{\sum_{P\in\mathcal{Q}^{j,\ell}}
  s_P\varphi_{2^{2j}I}\ast\psi_P}_{L^p}^q
    &\leq& C_N \left(\sum_{k\in\mathbb{Z}^d} \abs{s_{j,\ell,k}}^p \abs{Q_{2j}}^{-p}\abs{P_j}^{\frac{p}{2}}\abs{Q_{2j}}^{\frac{1}{2}}\right)^{q/p} \\
    &=& C_N \left(\sum_{k\in\mathbb{Z}^d} \abs{s_{j,\ell,k}}^p
        \abs{P_j}^{\frac{p}{2}}\abs{P_j}^{-\frac{2dp}{d+1}}\abs{P_j}^{\frac{d}{d+1}}\right)^{q/p}.
\end{eqnarray*}
If $1<p<\infty$ choose $N=a+b$ such that $a>d/p$ and $b>d(p-1)/p$.
Hölder's inequality gives
\begin{eqnarray*}
    && \!\!\!\!\!\!\!\!\!\!\!\!\!\!\!\!\!\!\!\!\!\!
        \left(\sum_{k\in\mathbb{Z}^d}\frac{\abs{s_{j,\ell,k}}\abs{Q_{2j}}^{-1}\abs{P_j}^{\frac{1}{2}}}
        {(1+2^j\abs{x-A^{-j}B^{-[\ell]}k})^N}\right)^p \\
    &\leq& \left(\sum_{k\in\mathbb{Z}^d}\frac{\abs{s_{j,\ell,k}}^p\abs{Q_{2j}}^{-p}\abs{P_j}^{\frac{p}{2}}}
        {(1+2^j\abs{x-A^{-j}B^{-[\ell]}k})^{ap}}\right)
        \left(\sum_{k\in\mathbb{Z}^d}\frac{1}
        {(1+2^j\abs{x-A^{-j}B^{-[\ell]}k})^{bp'}}\cdot\frac{2^{j bp'}}{2^{j bp'}}\right)^{p/p'} \\
    &\leq& C_{d,p} 2^{j bp} \left(\sum_{k\in\mathbb{Z}^d}\frac{\abs{s_{j,\ell,k}}^p\abs{Q_{2j}}^{-p}\abs{P_j}^{\frac{p}{2}}}
        {(1+2^j\abs{x-A^{-j}B^{-[\ell]}k})^{ap}}\right),
\end{eqnarray*}
because $\abs{A^{-j}B^{-[\ell]}x}\geq 2^{-2(j-1)}\abs{x}$ by Lemma
\ref{l:Nested_ellipsoids_HD}. Since $2^{j
bp}=\abs{P_j}^{-\frac{bp}{d+1}}$ and
$\abs{Q_{2j}}^{-p}=\abs{P_j}^{-\frac{2dp}{d+1}}$, we have
\begin{eqnarray*}
  \norm{\sum_{P\in\mathcal{Q}^{j,\ell}}s_P\varphi_{2^{2j}I}\ast\psi_P}_{L^p}^q
    &\leq& C_{d,p}\left(\sum_{k\in\mathbb{Z}^d}\abs{s_{j,\ell,k}}^p\abs{P_j}^{-\frac{bp}{d+1}-\frac{2dp}{d+1}+\frac{p}{2}}
        \abs{Q_{2j}}^{\frac{1}{2}}\right)^{q/p} \\
    &=& C_{d,p}\left(\sum_{k\in\mathbb{Z}^d}\abs{s_{j,\ell,k}}^p\abs{P_j}^{-\frac{bp}{d+1}-\frac{2dp}{d+1}+\frac{p}{2}+\frac{d}{d+1}}
        \right)^{q/p}.
\end{eqnarray*}
Let $\lambda=0$ if $0<p\leq 1$ or $\lambda>d(p-1)/p$ if
$1<p<\infty$. Then,
$$\norm{\sum_{P\in\mathcal{Q}^{j,\ell}}s_P\varphi_{2^{2j}I}\ast\psi_P}_{L^p}^q
\leq
C_{d,p}\left(\sum_{k\in\mathbb{Z}^d}[\abs{s_{j,\ell,k}}\abs{P_j}^{\frac{1}{2}-\frac{2d+\lambda}{d+1}+\frac{d}{p(d+1)}}]^p\right)^{q/p}.$$
Finally,
\begin{eqnarray*}
  \norm{f}_{\mathbf{B}_1}
    &\leq& C_{d,p} \left(\sum_{j=0}^\infty \abs{Q_{2j}}^{-\frac{\alpha_1
    q}{d}}\abs{Q_{2j}}^{-\frac{qs(d-1)}{d}}
        \sum_{\abs{[\ell]}\preceq 2^j}\norm{\sum_{P\in\mathcal{Q}^{j,\ell}} s_P\varphi_{2^{2j}I}\ast\psi_P}_{L^p}^q\right)^{1/q}  \\
    &\leq& C_{d,p} \left(\sum_{j=0}^\infty \abs{P_j}^{-\frac{2\alpha_1 q}{d+1}-\frac{2qs(d-1)}{d+1}}
        \sum_{\abs{[\ell]}\preceq 2^j}\left(\sum_{k\in\mathbb{Z}^d}
        \abs{s_{j,\ell,k}}^p\abs{P_j}^{p(\frac{1}{2}-\frac{2d+\lambda}{d+1}+\frac{d}{p(d+1)})}\right)^{q/p}\right)^{1/q}  \\
    &=& C_{d,p} \left(\sum_{j=0}^\infty\sum_{\abs{[\ell]}\preceq 2^j}
        \left(\sum_{k\in\mathbb{Z}^d}
        [\abs{s_{j,\ell,k}}\abs{P_j}^{\frac{1}{2}-\frac{2d+\lambda}{d+1}+\frac{d}{p(d+1)}-\frac{2(\alpha_1+s(d-1))}{d+1}}]^p\right)^{q/p}\right)^{1/q} \\
    &\leq& C_{d,p} \left(\sum_{j=0}^\infty\sum_{\abs{[\ell]}\preceq 2^j}
        \left(\sum_{k\in\mathbb{Z}^d}
        [\abs{s_{j,\ell,k}}\abs{P_j}^{-\alpha_2+\frac{d}{p(d+1)}-\frac{1}{2}}]^p\right)^{q/p}\right)^{1/q}
        =C_{d,p} \norm{\mathbf{s}}_{\mathbf{b}^{\alpha_2,q}_p},
\end{eqnarray*}
because $2d+\lambda+2(\alpha_1+s(d-1))<(d+1)(\alpha_2+1)$. Applying
Theorem \ref{t:Bndnss_T-S_Ops} finishes the proof.
\hfill $\blacksquare$ \vskip .5cm   

\vskip .5cm
\subsection{Vanishing norms of non-vanishing functions on Besov spaces}\label{sS:Non-vanishFncs_Bsv}
We now show that there exist sequences of non vanishing functions in
the norm of any of the shear anisotropic or isotropic spaces that
vanish in the norm of the other space.

\begin{thm}\label{t:Emmdngs_Bsv_in_BsvAB}
Let $\alpha_1,\alpha_2\in\mathbb{R}$ and
$0<p_1,p_2,q_1,q_2\leq\infty$ be such that
$2\alpha_1+\frac{3}{2}d<(d+1)(\alpha_2-\frac{d}{p_2(d+1)}+\frac{1}{2})+\frac{1}{2}+\frac{d}{p_1}$.
Then, there are sequences of functions in
$\mathbf{B}^{\alpha_2,q_2}_{p_2}(AB)$, with
$\norm{f^{(j)}}_{\mathbf{B}^{\alpha_2,q_2}_{p_2}(AB)}\approx 1$, for
all $j\in\mathbb{N}$, but $\lim_{j\rightarrow\infty}
\norm{f^{(j)}}_{\mathbf{B}^{\alpha_1,q_1}_{p_1}}\rightarrow 0$.
\end{thm}
\textbf{Proof}. For $j\geq 0$, let
$P_j=P_{j,0,0}\in\mathcal{Q}_{AB}$ and define
$\mathbf{s}^{(j)}=\{s^{(j)}_Q\}_{Q\in\mathcal{Q}_{AB}}$ such that
$$s_Q^{(j)}=\left\{\begin{array}{ccc}
                     0 & \text{if} & Q\neq P_j \\
                     \abs{P_j}^{\alpha_2-\frac{d}{p_2(d+1)}+\frac{1}{2}} & \text{if} &
                        Q=P_j.
                   \end{array}\right\}
$$
Then,
$\norm{\mathbf{s}^{(j)}}_{\mathbf{b}^{\alpha_2,q_2}_{p_2}(AB)}=1$,
for all $j\geq 0$. Thus, $f^{(j)}(x)=\sum_{Q\in\mathcal{Q}_{AB}}
s_Q^{(j)}\psi_Q(x)=\abs{P_j}^{\alpha_2-\frac{d}{p_2(d+1)}+\frac{1}{2}}\psi_{j,0,0}\in\mathbf{F}^{\alpha_2,q_2}_{p_2}(AB)$
with $\norm{f}_{\mathbf{F}^{\alpha_2,q_2}_{p_2}(AB)}\approx 1$, for
all $j\geq 0$. From the support conditions on $\hat{\varphi}$ and
$\hat{\psi}$, we have
\begin{eqnarray*}
  \norm{f^{(j)}}_{\mathbf{B}^{\alpha_1,q_1}_{p_1}}
    &=& \left(\sum_{\nu=(2j-5)_+}^{(2j-1)_+}(\abs{Q_\nu}^{-\frac{\alpha_1}{d}}
        \abs{P_j}^{\alpha_2-\frac{d}{p_2(d+1)}+\frac{1}{2}}
        \norm{\varphi_{2^\nu I}\ast
        \psi_{j,0,0}}_{L^{p_1}})^{q_1}\right)^{1/q_1}.
\end{eqnarray*}
Assume $\nu\sim 2j$. Since
$\varphi_{2^{2j}I}(x)=2^{2jd}\varphi(2^{2j}x)=\abs{Q_{2j}}^{-1}\varphi(2^{2j}x)$
and $\psi_{j,0,0}(x)=\abs{\text{det }A}^{\frac{j}{2}}\psi(A^j
x)=\abs{P_j}^{-\frac{1}{2}}\psi(A^j x)$, Lemma
\ref{l:conv_shrlts-wvlts} yields
$$\abs{\varphi_{2^{2j}I}\ast \psi_{j,0,0}(x)}
    \leq \frac{C_N \abs{Q_{2j}}^{-1}\abs{P_j}^{\frac{1}{2}}}{(1+2^j\abs{x})^N},$$
for some $C_N>0$ for all $N>d$. Taking $N>\max\{d,d/p_1\}$,
$\norm{\varphi_{2^{2j}I}\ast \psi_{j,0,0}(x)}_{L^{p_1}}\leq
C_{d,p_1}
\abs{Q_{2j}}^{-1}\abs{P_j}^{\frac{1}{2}}\abs{Q_{2j}}^{\frac{1}{2p_1}}$.
Hence,
\begin{eqnarray*}
  \norm{f^{(j)}}_{\mathbf{B}^{\alpha_1,q_1}_{p_1}}
    &\lesssim& C_{d,p_1} \abs{Q_{2j}}^{-\frac{\alpha_1}{d}}\abs{P_j}^{\alpha_2-\frac{d}{p_2(d+1)}+\frac{1}{2}}
        \cdot \abs{P_j}^{\frac{1}{2}}\abs{Q_{2j}}^{-1+\frac{1}{2p_1}}  \\
    &=&
    2^{-j(-2\alpha_1+(d+1)(\alpha_2-\frac{d}{p_2(d+1)}+\frac{1}{2})+\frac{d+1}{2}+2d(-1+\frac{1}{2p_1}))},
\end{eqnarray*}
tends to $0$, as $j\rightarrow\infty$, because
$2\alpha_1+\frac{3}{2}d<(d+1)(\alpha_2-\frac{d}{p_2(d+1)}+\frac{1}{2})+\frac{1}{2}+\frac{d}{p_1}$.

\hfill $\blacksquare$ \vskip .5cm   

\begin{thm}\label{t:Emmdngs_BsvAB_in_Bsv}
Let $\alpha_1,\alpha_2\in\mathbb{R}$ and
$0<p_1,p_2,q_1,q_2\leq\infty$ be such that $2\alpha_1+d
> (d+1)\alpha_2+\frac{d-1}{q_2}+\frac{2d}{p_1}$.
Then, there are sequences of functions in
$\mathbf{B}^{\alpha_1,q_1}_{p_1}$, with
$\norm{f^{(\nu)}}_{\mathbf{B}^{\alpha_1,q_1}_{p_1}}\approx 1$, for
all $\nu\in\mathbb{N}$, but $\lim_{\nu\rightarrow\infty}
\norm{f^{(\nu)}}_{\mathbf{B}^{\alpha_2,q_2}_{p_2}(AB)}\rightarrow
0$.
\end{thm}
\textbf{Proof}. For $\nu\geq 0$, let
$Q_\nu=Q_{\nu,0}\in\mathcal{D}_+$ and define
$\mathbf{s}^{(\nu)}=\{s^{(\nu)}_Q\}_{Q\in\mathcal{D}_+}$ such that
$$s_Q^{(\nu)}=\left\{\begin{array}{ccc}
                     0 & \text{if} & Q\neq Q_\nu \\
                     \abs{Q_\nu}^{\alpha_1-\frac{1}{p_1}+\frac{1}{2}} & \text{if} &
                        Q=Q_\nu.
                   \end{array}\right\}
$$
Then,
$\norm{\mathbf{s}^{(\nu)}}_{\mathbf{b}^{\alpha_1,q_1}_{p_1}}=1$, for
all $\nu\geq 0$. Thus, $f^{(\nu)}(x)=\sum_{Q\in\mathcal{D}_+}
s_Q^{(\nu)}\varphi_Q(x)=\abs{Q_\nu}^{\alpha_2-\frac{1}{p_2}+\frac{1}{2}}$
$\varphi_{\nu,0}\in\mathbf{F}^{\alpha_1,q_1}_{p_1}$ with
$\norm{f}_{\mathbf{F}^{\alpha_1,q_1}_{p_1}}\approx 1$, for all
$\nu\geq 0$. Consider the subsequence $f^{(2\nu)}$. From the compact
support conditions on $\hat{\varphi}$ and $\hat{\psi}$, we have
\begin{equation*}
    \norm{f^{(2\nu)}}_{\mathbf{B}_2}
    \lesssim \left(\sum_{\abs{[\ell]}\preceq 2^\nu} [\abs{P_\nu}^{-\alpha_2} \abs{Q_{2\nu}}^{\frac{\alpha_1}{d} -\frac{1}{p_1}+\frac{1}{2}}
        \norm{\varphi_{2\nu,0}\ast\psi_{A^{-\nu}B^{-\ell}}}_{L^{p_2}}]^{q_2}
        \right)^{1/{q_2}}.
\end{equation*}
Since
$\varphi_{2\nu,0}(x)=\abs{Q_{2\nu}}^{-\frac{1}{2}}\varphi(2^{2\nu}x)$
and $\psi_{A^{-\nu}B^{-\ell}}(x)=\abs{P_\nu}^{-1}\psi(B^\ell A^\nu
x)$, Lemma \ref{l:conv_shrlts-wvlts} yields
\begin{equation*}
    \abs{\varphi_{2\nu,0}\ast\psi_{A^{-\nu}B^{-\ell}}(x)}\leq
    \frac{C_N\abs{Q_{2\nu}}^{-\frac{1}{2}}}{(1+2^\nu\abs{x})^N},
\end{equation*}
for some $C_N>0$ for all $N>d$ and all $\ell$ such that
$\abs{[\ell]}\preceq 2^\nu$. Taking $N>\max\{d,d/p_2\}$,
$\norm{\varphi_{2\nu,0}\ast\psi_{A^{-\nu}B^{-\ell}}}_{L^{p_2}}\leq
C_{d,p_2}$. Therefore, since $\abs{\{\ell:\abs{[\ell]}\preceq
2^\nu\}}=2^{(\nu+1)(d-1)}+1\leq C_d 2^{\nu(d-1)}$, we finally get
\begin{eqnarray*}
  \norm{f^{(2\nu)}}_{\mathbf{B}_2}
    &\leq& C_{d,p_2} \left(\sum_{\abs{[\ell]}\preceq 2^\nu}
        [\abs{P_\nu}^{-\alpha_2}\abs{Q_{2\nu}}^{\frac{\alpha_1}{d}-\frac{1}{p_1}+\frac{1}{2}}]^{q_2}\right)^{1/q_2} \\
    &\leq& C_{d,p_2} \left([2^{\frac{\nu(d-1)}{q_2}}2^{\nu(d+1)\alpha_2}2^{-2\nu d(\frac{\alpha_1}{d}-\frac{1}{p_1}+\frac{1}{2})}]^{q_2}\right)^{1/q_2}  \\
    &\leq& C_{d,p_2}
        2^{-\nu[-\frac{d-1}{q_2}-(d+1)\alpha_2+2d(\frac{\alpha_1}{d}-\frac{1}{p_1}+\frac{1}{2})]},
\end{eqnarray*}
which tends to $0$, as $\nu\rightarrow\infty$, if $2\alpha_1+d
> (d+1)\alpha_2+\frac{d-1}{q_2}+\frac{2d}{p_1}$.
\hfill $\blacksquare$ \vskip .5cm   

\vskip 0.5cm
\subsection{Embeddings of Triebel-Lizorkin spaces}\label{sS:Emmbeddings_T-L}
\begin{thm}\label{t:ClssDydc-T-L_in_ShrAnistr-T-L}
Let $\alpha_1,\alpha_2\in \mathbb{R}$, $0<q\leq\infty$, $0<p<\infty$
and $\lambda>d\max(1,1/q,1/p)$. If
$(d+1)\alpha_2+\frac{d-1}{q}+\lambda\leq 2\alpha_1$,
$$\mathbf{F}^{\alpha_1,q}_p \hookrightarrow \mathbf{F}^{\alpha_2,q}_p(AB).$$
\end{thm}
\textbf{Proof}. We only prove for $q<\infty$, case $q=\infty$ is
similar. To shorten notation write
$\mathbf{F}_1=\mathbf{F}^{\alpha_1,q}_p$,
$\mathbf{f}_1=\mathbf{f}^{\alpha_1,q}_p$ and
$\mathbf{F}_2=\mathbf{F}^{\alpha_2,q}_p(AB)$. From the compact
support conditions on $\hat{\varphi}$ and $\hat{\psi}$, and since
$\varphi_{2j,k}(x)=\abs{Q_{2j}}^{-\frac{1}{2}}\varphi(2^{2j}x-k)$
and $\psi_{A^{-j}B^{-\ell}}(x)=\abs{P_j}^{-1}\psi(B^\ell A^j x)$,
Lemma \ref{l:conv_shrlts-wvlts} yields
\begin{equation*}
    \abs{\varphi_{2j,k}\ast\psi_{A^{-j}B^{-\ell}}(x)}\leq \frac{C_N
    \abs{Q_{2j}}^{-\frac{1}{2}}}{(1+2^j\abs{x+2^{-2j}k})^N},
\end{equation*}
for some $C_N>0$ for all $N>d$ and all $\ell$ such that
$\abs{[\ell]}\preceq 2^j$. Then,
\begin{eqnarray*}
  \norm{f}_{\mathbf{F}_2}
    &\lesssim& C_N \norm{\left(\sum_{j\geq 0}\abs{P_j}^{-\alpha_2q}\sum_{\abs{[\ell]}\preceq 2^j}
        [\sum_{k\in\mathbb{Z}^d}
        \abs{s_{2j,k}}\frac{\abs{Q_{2j}}^{-\frac{1}{2}}}{(1+2^j\abs{\cdot+2^{-2j}k})^N}]^q\right)^{1/q}}_{L^p},
\end{eqnarray*}
for all $N>d$. Let $\lambda>d\max(1,1/q,1/p)$. Following the proof
of the second part of Theorem \ref{t:S-T-psi-ops_Bnd}, if $x\in Q$
and $Q\in\mathcal{D}^{2j}$,
\begin{eqnarray*}
    & &\!\!\!\!\!\!\!\!\!\!\!\!\!\!\!\!\!\!\!\!\!\!\!\!\!\!\!\!\!\!\!\!
        [\sum_{k\in\mathbb{Z}^d}
        \frac{\abs{s_{2j,k}}\abs{Q_{2j}}^{-\frac{1}{2}}}{(1+2^j\abs{x-2^{-2j}k})^\lambda}]^q\\
    &=& [\sum_{k\in\mathbb{Z}^d} \frac{\abs{s_{2j,k}}\abs{Q_{2j}}^{-\frac{1}{2}}
        \cdot
        2^{j\lambda}}{2^{j\lambda}(1+2^j\abs{x-2^{-2j}k})^\lambda}]^q
        \leq [2^{j\lambda}\sum_{k\in\mathbb{Z}^d}
        \frac{\abs{s_{2j,k}}\abs{Q_{2j}}^{-\frac{1}{2}}}{(1+2^{2j}\abs{x-2^{-2j}k})^\lambda}]^q\\
    &\lesssim& [2^{j\lambda}\sum_{Q\in\mathcal{D}^{2j}}
    \abs{Q}^{-\frac{1}{2}}\abs{(s_{1,\lambda}^\ast)_Q}\chi_Q(x)]^q
        = 2^{j\lambda q}\sum_{Q\in\mathcal{D}^{2j}}
        [\abs{(s_{1,\lambda}^\ast)_Q}\tilde{\chi}_Q(x)]^q,
\end{eqnarray*}
because $\mathcal{D}^{2j}$ is a partition of $\mathbb{R}^d$. Since
$\abs{\{\ell:\abs{[\ell]}\preceq 2^j\}}=2^{(j+1)(d-1)}+1\leq C_d
2^{j(d-1)}$, we have
\begin{eqnarray*}
  \norm{f}_{\mathbf{F}_2}
    &\lesssim& C_{d,p,q}\norm{\left(\sum_{j\geq 0}\abs{P_j}^{-\alpha_2q}2^{j(d-1)}2^{j\lambda q}
        \sum_{Q\in\mathcal{D}^{2j}}
        [\abs{(s_{1,\lambda}^\ast)_Q}\tilde{\chi}_Q(\cdot)]^q\right)^{1/q}}_{L^p} \\
    &\leq& C_{d,p,q} \norm{\left(\sum_{j\geq 0}\sum_{Q\in\mathcal{D}^{2j}}
        [\abs{Q}^{-\frac{\alpha_1}{d}}\abs{(s_{1,\lambda}^\ast)_Q}\tilde{\chi}_Q(\cdot)]^q\right)^{1/q}}_{L^p}\\
    &\leq& C_{d,p,q} \norm{\left(\sum_{j\geq 0}\sum_{Q\in\mathcal{D}^j}
        [\abs{Q}^{-\frac{\alpha_1}{d}}\abs{(s_{1,\lambda}^\ast)_Q}\tilde{\chi}_Q(\cdot)]^q\right)^{1/q}}_{L^p}
        = \norm{\mathbf{s}^\ast_{1,\lambda}}_{\mathbf{f}_1},
\end{eqnarray*}
because $(d+1)j\alpha_2+\frac{(d-1)j}{q}+j\lambda \leq 2j\alpha_1$.
Following the proof of Lemma 2.3 of \cite{FJ90} (where the
restriction on $\lambda$ is used) it can be concluded that
$\norm{\mathbf{s}_{1,\lambda}^\ast}_{\mathbf{f}_1}
    \lesssim\norm{\mathbf{s}}_{\mathbf{f}_1}$,
and from Theorem 2.2 in \cite{FJ90},
$\norm{\mathbf{s}}_{\mathbf{f}_1}\lesssim\norm{f}_{\mathbf{F}_1}$,
which finishes the proof.
\hfill $\blacksquare$ \vskip .5cm   

\begin{thm}\label{t:ShrAnistr-T-L_in_ClssDydc-T-L}
Let $\alpha_1,\alpha_2\in \mathbb{R}$, $0<q\leq\infty$ and
$0<p<\infty$ be such that $2\alpha_1+d+(d-1)(1-1/q)_+\leq
(d+1)\alpha_2+1$, where $(1-1/q)_+=\max\{0,1-1/q\}$. Then,
$$\mathbf{F}^{\alpha_2,q}_p(AB) \hookrightarrow \mathbf{F}^{\alpha_1,q}_p.$$
\end{thm}
\textbf{Proof.} We only prove for $q<\infty$, case $q=\infty$ is
similar. To shorten notation write
$\mathbf{F}_1=\mathbf{F}^{\alpha_1,q}_p$,
$\mathbf{F}_2=\mathbf{F}^{\alpha_2,q}_p(AB)$ and
$\mathbf{f}_2=\mathbf{f}^{\alpha_2,q}_p(AB)$. Suppose
$f=\sum_{P\in\mathcal{Q}_{AB}}s_P\psi_P\in \textbf{F}_2$. From the
compact support conditions on $(\varphi_{2^\nu I})^\wedge$ and
$(\psi_{j,\ell,k})^\wedge$, we have
\begin{eqnarray*}
  \norm{f}_{\textbf{F}_1}
    &=& \norm{\left(\sum_{\nu\geq 0} (\abs{Q_\nu}^{-\frac{\alpha_1}{d}}\abs{\varphi_{2^\nu I}\ast f(\cdot)})^q\right)^{1/q}}_{L^p} \\
    &\lesssim& \norm{\left(\sum_{\nu\geq 0} \abs{Q_{2\nu}}^{-\frac{\alpha_1 q}{d}}
        (\sum_{\abs{[\ell]}\preceq 2^\nu}\sum_{k\in\mathbb{Z}^d}\abs{s_{\nu,\ell,k}}
        \abs{\varphi_{2^{2\nu}I}\ast \psi_{\nu,\ell,k}(\cdot)})^q\right)^{1/q}}_{L^p} \\
    &\lesssim& C_N \norm{\left(\sum_{\nu\geq 0} \abs{Q_{2\nu}}^{-\frac{\alpha_1 q}{d}}
        (\sum_{\abs{[\ell]}\preceq2^\nu}\sum_{k\in\mathbb{Z}^d}\abs{s_{\nu,\ell,k}}
        \frac{\abs{Q_{2\nu}}^{-1} \abs{P_\nu}^{-1/2} \abs{P_\nu}}{(1+2^\nu\abs{\cdot-A^{-\nu}B^{-\ell}k})^N})^q\right)^{1/q}}_{L^p},
\end{eqnarray*}
for some $C_N>0$ for all $N>d$, by Lemma \ref{l:conv_shrlts-wvlts}.
Continuing as in the second part of the proof of Theorem
\ref{t:S-T-psi-ops_Bnd}, if $x\in P$ and
$P\in\mathcal{Q}^{\nu,\ell}$,
\begin{eqnarray*}
  \norm{f}_{\textbf{F}_1}
    &\lesssim& C_N \norm{\left(\sum_{\nu\geq 0} \abs{Q_{2\nu}}^{-\frac{\alpha_1 q}{d}-q} \abs{P_\nu}^q
        [\sum_{\abs{[\ell]}\preceq2^\nu}\sum_{P\in\mathcal{Q}^{\nu,\ell}}\abs{s_{P}}\cdot
        \frac{\abs{P}^{-1/2}}{(1+2^\nu\abs{\cdot-x_P})^N}]^q\right)^{1/q}}_{L^p} \\
    &\lesssim& C_N \norm{\left(\sum_{\nu\geq 0} \abs{P_\nu}^{\frac{2dq}{d+1}(-\frac{\alpha_1 }{d}-1)+q}
        [\sum_{\abs{[\ell]}\preceq2^\nu}\sum_{P\in\mathcal{Q}^{\nu,\ell}}
        (s_{1,N}^\ast)_P\tilde{\chi}_P(\cdot)]^q\right)^{1/q}}_{L^p},
\end{eqnarray*}
because $\mathcal{Q}^{\nu,\ell}$ is a partition of $\mathbb{R}^d$.
However,
$\sum_{\abs{[\ell]}\preceq2^\nu}\sum_{P\in\mathcal{Q}^{\nu,\ell}}\chi_P$
is not a partition of $\mathbb{R}^d$. If $0<q\leq 1$ we use the
$q$-triangle inequality $\abs{a+b}^q\leq \abs{a}^q+\abs{b}^q$
($N>d/q$) to get
\begin{equation*}
    [\sum_{\abs{[\ell]}\preceq 2^\nu}\sum_{P\in\mathcal{Q}^{\nu,\ell}} (s^\ast_{1,N})_P\tilde{\chi}_P(\cdot)]^q
    \leq  \sum_{\abs{[\ell]}\preceq 2^\nu}\sum_{P\in\mathcal{Q}^{\nu,\ell}}
    [(s^\ast_{1,N})_P\tilde{\chi}_P(\cdot)]^q,
\end{equation*}
or Hölder's inequality if $1<q$ ($N>d$) to get
\begin{eqnarray*}
  [\sum_{\abs{[\ell]}\preceq 2^\nu}\sum_{P\in\mathcal{Q}^{\nu,\ell}} (s^\ast_{1,N})_P\tilde{\chi}_P(\cdot)]^q
    &\leq& C_d 2^{\nu(d-1)q(1-\frac{1}{q})} \sum_{\abs{[\ell]}\preceq 2^\nu}(\sum_{P\in\mathcal{Q}^{\nu,\ell}}
        (s^\ast_{1,N})_P\tilde{\chi}_P(\cdot))^q \\
    &=& C_d 2^{\nu(d-1)q(1-\frac{1}{q})} \sum_{\abs{[\ell]}\preceq 2^\nu}\sum_{P\in\mathcal{Q}^{\nu,\ell}}
        [(s^\ast_{1,N})_P\tilde{\chi}_P(\cdot)]^q,
\end{eqnarray*}
because $\mathcal{Q}^{\nu,\ell}$ is a partition of $\mathbb{R}^d$.
Let $\lambda>(d+1)\max(1,1/q,1/p)$,
\begin{eqnarray*}
  \norm{f}_{\textbf{F}_1}
    &\lesssim& \norm{\left(\sum_{P\in\mathcal{Q}_{AB}}
        [\abs{P_\nu}^{\frac{2d}{d+1}(-\frac{\alpha_1}{d}-1)+1} 2^{\nu(d-1)(1-\frac{1}{q})_+} (s_{1,\lambda}^\ast)_P\tilde{\chi}_P(\cdot)]^q\right)^{1/q}}_{L^p} \\
    &\lesssim& \norm{\left(\sum_{P\in\mathcal{Q}_{AB}}
        [\abs{P}^{-\alpha_2}(s_{1,\lambda}^\ast)_P\tilde{\chi}_P(\cdot)]^q\right)^{1/q}}_{L^p}
        =\norm{\mathbf{s}_{1,\lambda}^\ast}_{\textbf{f}_2}.
\end{eqnarray*}
By Lemma \ref{l:s_ast_bnd_s} and Theorem \ref{t:S-T-psi-ops_Bnd} the
proof is complete.

\hfill $\blacksquare$ \vskip .5cm   

\subsection{Vanishing norms of non vanishing functions on Triebel-Lizorkin spaces}\label{sS:More_Relations_T-L-Shrl_T-L_Dyadic}
We now show that there exist sequences of non vanishing functions in
the norm of any of the shear anisotropic or isotropic spaces that
vanish in the norm of the other space.

\begin{thm}\label{t:Dydc-T-L_fade_with_Shrlt-T-L}
Let $\alpha_1,\alpha_2\in \mathbb{R}$, $0<q_1,q_2\leq\infty$ and
$0<p_1,p_2<\infty$. Then, there exist sequences of functions
$\{f^{(j)}\}_{j\geq 0}$ such that
$\norm{f^{(j)}}_{\mathbf{F}^{\alpha_2,q_2}_{p_2}(AB)}\approx 1$, but
that $\norm{f^{(j)}}_{\mathbf{F}^{\alpha_1,q_1}_{p_1}}\rightarrow
0$, $j\rightarrow\infty$, if
$2(\alpha_1+d)<(d+1)(\alpha_2-\frac{1}{p_2}+1)+\frac{d}{p_1}$.
\end{thm}
\textbf{Proof}. For $j\geq 0$, let
$P_j=P_{j,0,0}\in\mathcal{Q}_{AB}$ and define
$\mathbf{s}^{(j)}=\{s_Q^{(j)}\}_{Q\in\mathcal{Q}_{AB}}$ such that
$$s_Q^{(j)}=\left\{\begin{array}{ccc}
                     0 & \text{if} & Q\neq P_j \\
                     \abs{P_j}^{\alpha_2-\frac{1}{p_2}+\frac{1}{2}} & \text{if} & Q=
                     P_j.
                   \end{array}
\right\}$$ Then,
$\norm{\mathbf{s}^{(j)}}_{\mathbf{f}^{\alpha_2,q_2}_{p_2}}=1$, for
all $j\geq 0$. Thus,
$f^{(j)}(x)=\sum_{Q\in\mathcal{Q}_{AB}}s_Q^{(j)}\psi_Q(x)
    =\abs{P_j}^{\alpha_2-\frac{1}{p_2}+\frac{1}{2}}$ $\psi_{j,0,0}(x)\in\mathbf{F}^{\alpha_2,q_2}_{p_2}(AB)$
with $\norm{f}_{\mathbf{F}^{\alpha_2,q_2}_{p_2}(AB)}\approx 1$. From
the compact support conditions on $\hat{\psi}$ and $\hat{\varphi}$,
we have
\begin{eqnarray*}
  \norm{f^{(j)}}_{\mathbf{F}^{\alpha_1,q_1}_{p_1}}
    &=& \norm{\left(\sum_{\nu=0}^\infty [2^{\nu\alpha_1}
        \abs{f^{(j)}\ast\varphi_{2^\nu I}}]^{q_1}\right)^{1/q_1}}_{L^{p_1}} \\
    &\lesssim& \norm{2^{2j\alpha_1}
        \abs{f^{(j)}\ast\varphi_{2^{2j}
        I}}}_{L^{p_1}}.
\end{eqnarray*}
Lemma \ref{l:conv_shrlts-wvlts} yields
\begin{eqnarray*}
  \abs{f^{(j)}\ast\varphi_{2^{2j}I}(x)}
    &=& \abs{\int_{\mathbb{R}^d} \abs{P_j}^{\alpha_2-\frac{1}{p_2}+\frac{1}{2}} \abs{\text{det }A}^{j/2}
        \psi(A^j(x-y))2^{2jd}\varphi(2^{2j}y)dy} \\
    &\lesssim& \abs{P_j}^{(\alpha_2-\frac{1}{p_2})} 2^{2jd}
        \int_{\mathbb{R}^d} \abs{\psi(A^j(x-y))}\abs{\varphi(2^{2j}y)}
        dy\\
    &\lesssim& \frac{\abs{P_j}^{(\alpha_2-\frac{1}{p_2}+1)} 2^{2jd}}{(1+2^j\abs{x})^N},
\end{eqnarray*}
for every $N>d$. Then, for $N$ such that $Np_1>d$, we have
\begin{eqnarray*}
  \norm{f^{(j)}}_{\mathbf{F}^{\alpha_1,q_1}_{p_1}}
    &\leq& C_{N,q_1} \left(\int_{\mathbb{R}^d} 2^{2j\alpha_1p_1}\cdot
        \frac{[\abs{P_j}^{(\alpha_2-\frac{1}{p_2}+1)} 2^{2jd}]^{p_1}}{(1+2^j\abs{x})^{Np_1}} dx\right)^{1/p_1} \\
    &=& C_{N,p_1,q_1}
        2^{2j\alpha_1-(d+1)j(\alpha_2-\frac{1}{p_2}+1)+2jd-\frac{dj}{p_1}},
\end{eqnarray*}
which tends to $0$ as $j\rightarrow\infty$ if
$2(\alpha_1+d)<(d+1)(\alpha_2-\frac{1}{p_2}+1)+\frac{d}{p_1}$.
\hfill $\blacksquare$ \vskip .5cm   

\begin{thm}\label{t:ShrAnistr-T-L_in_ClssDydc-T-L}
Let $\alpha_1,\alpha_2\in \mathbb{R}$, $0<q_1,q_2\leq\infty$ and
$0<p_1,p_2<\infty$. Then, there exist sequences of functions
$\{f^{(\nu)}\}_{\nu\geq 0}$ such that
$\norm{f^{(\nu)}}_{\mathbf{F}^{\alpha_1,q_1}_{p_1}}\approx 1$, but
that $\norm{f^{(\nu)}}_{\mathbf{F}^{\alpha_2,q_2}_{p_2}(AB)}$
$\rightarrow 0$, $\nu\rightarrow\infty$, if
$\frac{d-1}{q_2}+(d+1)\alpha_2+\frac{2d}{p_1}<2\alpha_1+\frac{d}{p_2}$.
\end{thm}
\textbf{Proof}. For a sequence
$\mathbf{s}^{(\nu)}=\{s_{\nu,0}\}_{j\geq 0}$ such that
$\abs{s_{\nu,0}}=\abs{Q_\nu}^{\frac{\alpha_1}{d}-\frac{1}{p_1}+\frac{1}{2}}$,
$\norm{\mathbf{s}^{(\nu)}}_{\mathbf{f}^{\alpha_1,q_1}_{p_1}}=1$, for
all $\nu\geq 0$. This means that
$f^{(\nu)}(x)=s_{\nu,0}\varphi_{\nu,0}(x)\in
\mathbf{F}^{\alpha_1,q_1}_{p_1}$ and
$\norm{f^{(\nu)}}_{\mathbf{F}^{\alpha_1,q_1}_{p_1}}\approx 1$.
Consider the subsequence $f^{(2j)}$. The conditions on the compact
support of $\hat{\varphi}$ and $\hat{\psi}$ give
\begin{equation*}
    \norm{f^{(2j)}}_{\mathbf{F}_2}
    \lesssim \norm{\left(\sum_{\abs{[\ell]}\preceq 2^j}
        [\abs{P_j}^{-\alpha_2}\abs{Q_{2j}}^{\frac{\alpha_1}{d}-\frac{1}{p_1}+\frac{1}{2}}
            \abs{\varphi_{2j,0}\ast\psi_{A^{-j}B^{-\ell}}}]^{q_2}\right)^{1/q_2}}_{L^{p_2}}.
\end{equation*}
Lemma \ref{l:conv_shrlts-wvlts} yields
\begin{eqnarray*}
  \abs{\varphi_{2j,0}\ast\psi_{A^{-j}B^{-\ell}}(x)}
    &=& \abs{\int_{\mathbb{R}^2} 2^{2jd/2}\varphi(2^{2j}y)\abs{\text{det }A}^j\psi(B^\ell A^j(x-y)) dy} \\
    &\lesssim& \abs{Q_{2j}}^{-\frac{1}{2}} \abs{P_j}^{-1}
        \int_{\mathbb{R}^2} \abs{\psi(B^\ell A^j(x-y))}\abs{\varphi(2^{2j}y)}dy \\
    &\lesssim&
       \frac{\abs{Q_{2j}}^{-\frac{1}{2}}}{(1+2^j\abs{x})^N},
\end{eqnarray*}
for every $N>d$ and all $\ell$ such that $\abs{[\ell]}\preceq 2^j$.
Then, since $\abs{\{\ell:\abs{[\ell]}\preceq 2^j\}}\lesssim
c_d2^{j(d-1)}$ for $N>d$ such that $Np_2>d$,
\begin{eqnarray*}
  \norm{f^{(2j)}}_{\mathbf{F}_2}
   &\leq& 2^{\frac{j(d-1)}{q_2}}\abs{P_j}^{-\alpha_2} \abs{Q_{2j}}^{\frac{\alpha_1}{d}-\frac{1}{p_1}}
        \left(\int_{\mathbb{R}^2} \frac{dx}{(1+2^j\abs{x})^{Np_2}} \right)^{1/p_2} \\
    &\lesssim& 2^{\frac{j(d-1)}{q_2}}2^{j(d+1)\alpha_2-2j\alpha_1+\frac{2jd}{p_1}-\frac{jd}{p_2}},
\end{eqnarray*}
which tends to $0$, as $j\rightarrow\infty$, if
$\frac{d-1}{q_2}+(d+1)\alpha_2+\frac{2d}{p_1}<2\alpha_1+\frac{d}{p_2}$.
\hfill $\blacksquare$ \vskip .5cm   

\vskip 1cm
\section{Proofs}\label{S:Proofs}

\vskip 0.5cm
\subsection{Proofs for Section \ref{S:Basic_Results}}\label{sS:Proofs_Basic_Results}
In order to prove Lemma \ref{l:c:l:conv_shrlts} we need two previous
results.

\begin{lem}\label{l:Nested_ellipsoids_HD}
Let $A^j$ and $B^{[\ell]}$ be as in Section \ref{S:Intro}. Then,
$$C_d2^j\abs{x}< \abs{B^{[\ell]} A^j x},$$
for every $j\geq 0$, $\abs{[\ell]}\preceq 2^j$ and all
$x\in\mathbb{R}^d$ with $C_d=2^{-d+1}$.
\end{lem}
\noindent\textbf{Proof}. It suffices to prove the lemma for
$x\in\partial\mathbb{U}:=\{y\in\mathbb{R}^d:y_1^2+\ldots+y_d^2=1\}$.
From (\ref{e:BAx}) and since
$\norm{\cdot}_{\ell^1(\mathbb{R}^{d-1})}\leq
\sqrt{d-1}\norm{\cdot}_{\ell^2(\mathbb{R}^{d-1})}$ and
$\abs{\ell_n}\leq 2^j$, $n=1,\ldots,d-1$, we have
\begin{eqnarray*}
  \abs{B^{[\ell]} A^j x}
    &\geq& \abs{2^{2j}x_1+2^j\ell_1x_2+\ldots+2^j\ell_{d-1}x_d} \\
    &\geq& \abs{2^{2j}\abs{x_1}-\abs{2^j\ell_1x_2+\ldots+2^j\ell_{d-1}x_d}} \\
    &\geq& 2^{2j}\abs{\abs{x_1}-\sqrt{d-1}(\sum_{n=2}^d\abs{x_n}^2)^{1/2}} \\
    &=& 2^{2j}\abs{\abs{x_1}-\sqrt{d-1}(1-\abs{x_1}^2)^{1/2}},
\end{eqnarray*}
because $x_1^2+\ldots+x_d^2=1$. Consider $\abs{x_1}^2\geq
(2^{2(d-1)}-1)/2^{2(d-1)}$. Then,
\begin{eqnarray*}
  \abs{B^{[\ell]} A^j x}
    &\geq& 2^{2j} \left(\sqrt{\frac{2^{2(d-1)}-1}{2^{2(d-1)}}}-\sqrt{d-1}\sqrt{1-\frac{2^{2(d-1)}-1}{2^{2(d-1)}}}\right) \\
    &\geq& 2^{2j} \left(\frac{2^{2(d-1)}-d}{2^{2(d-1)}+2^{(d-1)}\sqrt{d}}\right)
        \geq 2^{2(j-1)},
\end{eqnarray*}
because $d\geq 2$. When $\abs{x_1}^2< (2^{2(d-1)}-1)/2^{2(d-1)}$,
$x\in\partial \mathbb{U}$ implies
$\abs{x_2}^2+\cdots+\abs{x_d}^2>2^{-2(d-1)}$. Therefore,
$\abs{B^{[\ell]}A^jx}>2^{j-(d-1)}$.

Similarly one can prove $A^{-j}B^{-[\ell]} x> 2^{-2(j-1)}\abs{x}$,
$j\geq 0$, $[\ell]\preceq 2^j$.
\hfill $\blacksquare$ \vskip .5cm   

\begin{lem}\label{l:conv_shrlts}
Let $g,h\in\mathcal{S}$. Then, for every $N>d$, $i=j-1, j, j+1\geq
0$, $\abs{[m]}\preceq 2^i$ and $\abs{[\ell]}\preceq 2^j$ there exist
$C_N>0$ such that
$$\abs{g_{j,\ell,k}\ast h_{i,m,n}(x)}\leq
\frac{C_N}{(1+2^i\abs{x-A^{-i}B^{-m}n-A^{-j}B^{-\ell}k})^{N}},$$ for
all $x\in\mathbb{R}^d$.
\end{lem}
\noindent\textbf{Proof.} Since $g,h\in\mathcal{S}$, there exists
$C_N>0$ such that $\abs{g(x)}, \abs{h(x)}\leq
\frac{C_N}{(1+\abs{x})^N}$ for all $N\in\mathbb{N}$. Then,
\begin{eqnarray*}
  \abs{g_{j,\ell,k}\ast h_{i,m,n}(x)}
   &\leq& \abs{\text{ det }A}^{(j+i)/2} \int_{\mathbb{R}^d} \frac{C_N}{(1+\abs{B^\ell A^jy})^N}\cdot \frac{C_N}{(1+\abs{B^m A^i(x'-y)})^N} dy
\end{eqnarray*}
where $x'=x-A^{-j}B^{-\ell}k-A^{-i}B^{-m} n$. Notice that, since
$i=j-1,j,j+1$, $\abs{\text{ det }A}^{(j+i)/2}\simeq \abs{\text{ det
}A}^{j}\simeq \abs{\text{ det }A}^{i}$. Following \cite[\S 6]{HW96},
define
\begin{eqnarray*}
  E_1 &=& \{y\in\mathbb{R}^d:\abs{B^m A^i(x'-y)}\leq 3\} \\
  E_2 &=& \{y\in\mathbb{R}^d:\abs{B^m A^i(x'-y)}> 3, \abs{y}\leq \abs{x'}/2\} \\
  E_3 &=& \{y\in\mathbb{R}^d:\abs{B^m A^i(x'-y)}> 3, \abs{y}> \abs{x'}/2\}.
\end{eqnarray*}
Lemma \ref{l:Nested_ellipsoids_HD} yields the next three bounds. For
$y\in E_1$ we have $1+2^i\abs{x'}\leq 1+C_d^{-1}\abs{B^m
A^i(x'-y)}+2^i\abs{y}\leq 1+3C_d^{-1}+2^{j+1}\abs{y}\leq
1+3C_d^{-1}+2C_d^{-1}\abs{B^\ell A^j y} \leq c_d(1+\abs{B^\ell A^j
y})$. If $y\in E_3$, $1+2^i\abs{x'}\leq 1+2^{j+2}\abs{y}\leq
1+2^2C_d^{-1}\abs{B^\ell A^jy}
 \leq c_d(1+\abs{B^\ell A^jy})$. When $y\in E_2$, $2^{i-1}\abs{x'}\leq
2^i\abs{x'}-2^i\abs{y}<2^i\abs{x'-y}$, which implies $4\abs{B^m
A^i(x'-y)}=\abs{B^m A^i(x'-y)}+3\abs{B^m A^i(x'-y)}\geq
3+3C_d2^i\abs{x'-y}\geq c'_d(1 +2^i\abs{x'-y})\geq
c_d(1+2^i\abs{x'})$. Thus, since $\abs{j-i}\leq 1$,
\begin{eqnarray*}
  \abs{g_{j,\ell,k}\ast h_{i,m,n}(x)}
    &\lesssim& \frac{C_N \abs{\text{ det }A}^{i}}{(1+2^i\abs{x'})^N}\int_{E_1\cup E_3}\frac{C_N}{(1+\abs{B^m A^i(x'-y)})^N} dy \\
    & & \;\;\;\; + \frac{C_N \abs{\text{ det }A}^{i}}{(1+2^i\abs{x'})^N}\int_{E_2}\frac{C_N}{(1+\abs{B^\ell A^j y})^N} dy\\
    &\lesssim& \frac{C_N}{(1+2^i\abs{x'})^N}
\end{eqnarray*}
for some $C_N>0$ for every $N>d$, doing a change of variables to
bound the integrals with a constant independent of $i, j , l$ and
$m$. The result follows by replacing back
$x'=x-A^{-j}B^{-\ell}k-A^{-i}B^{-m} n$ in the estimate above.

\hfill $\blacksquare$ \vskip .5cm   

\noindent As a corollary for Lemma \ref{l:conv_shrlts} we have our
first ``almost orthogonality" property for the anisotropic and shear
dilations for functions in $\mathcal{S}$.

\noindent\textbf{Proof of Lemma \ref{l:c:l:conv_shrlts}.} Identify
$(j,\ell,0)$ with $P$ and $(i,m,n)$ with $Q$. Write
$\abs{g_{A^{-j}B^{-\ell}}\ast
h_{i,m,n}(x)}=\abs{\abs{P}^{-1/2}g_P\ast h_Q}$. Since $\abs{i-j}\leq
1$, $\abs{P}^{-1/2}\thicksim\abs{Q}^{-1/2}$. Then, Lemma
\ref{l:conv_shrlts} yields
$$\abs{\abs{P}^{-1/2}g_P\ast h_Q}\leq \frac{C_N \abs{P}^{-1/2}}
    {(1+2^j\abs{x-x_Q})^N}\lesssim \frac{C_N \abs{Q}^{-1/2}}{(1+2^j\abs{x-x_Q})^N}.$$

\hfill $\blacksquare$ \vskip .5cm   

\vskip .5cm We now present the proof for the second ``almost
orthogonality" result regarding dyadic isotropic dilated function
and shear anisotropic dilated function.

\noindent\textbf{Proof of Lemma \ref{l:conv_shrlts-wvlts}}. Since
$\psi,\varphi\in\mathcal{S}$,
\begin{eqnarray*}
    & & \!\!\!\!\!\!\!\!\!\!\!\!\!\!\!\!\!\!\!\!\!\!\!\!\!\!\!\!\!\!\!\!
        \int_{\mathbb{R}^d} \abs{\psi(B^\ell
        A^j(x-y))}\abs{\varphi(2^{2j}y)} dy\\
    &\lesssim& \int_{\mathbb{R}^d}
        \frac{1}{(1+\abs{B^\ell A^j(x-y)})^N}
        \frac{1}{(1+\abs{2^{2j}y})^N}dy.
\end{eqnarray*}
Define
\begin{eqnarray*}
  E_1 &=& \{y\in\mathbb{R}^d:\abs{y}> \frac{\abs{x}}{2}\} \\
  E_2 &=& \{y\in\mathbb{R}^d:\abs{y}\leq \frac{\abs{x}}{2}\} \\
\end{eqnarray*}
If $y\in E_1$, $1+2^j\abs{x}\leq 1+2^{j+1}\abs{y}\leq
2(1+2^{2j}\abs{y})$. When $y\in E_2$, $\frac{1}{2}\abs{x}<
\abs{x}-\abs{y}\leq \abs{x-y}$, which implies $4(1+\abs{B^\ell
A^j(x-y)})\geq 1+4\abs{B^\ell A^j (x-y)}\geq 1+4C_d 2^j\abs{x-y}\geq
c_d(1+2^j\abs{x})$, by Lemma \ref{l:Nested_ellipsoids_HD}. Hence,
\begin{eqnarray*}
    & & \!\!\!\!\!\!\!\!\!\!\!\!\!\!\!\!\!\!\!\!\!\!\!\!\!\!\!\!\!\!\!\!
        \int_{\mathbb{R}^d} \abs{\psi(B^\ell
        A^j(x-y))}\abs{\varphi(2^{2j}y)} dy\\
    &\lesssim& \frac{1}{(1+2^j\abs{x})^N} \int_{E_1}
        \frac{1}{(1+\abs{B^\ell A^j(x-y)})^N}dy \\
    & & \;\;\; + \frac{1}{(1+2^j\abs{x})^N} \int_{E_2}
        \frac{1}{(1+2^{2j}\abs{y})^N}dy\\
    &\lesssim& \left[\frac{2^{-(d+1)j}}{(1+2^j\abs{x})^N}+\frac{2^{-2dj}}{(1+2^j\abs{x})^N}\right]
        \lesssim \frac{2^{-(d+1)j}}{(1+2^j\abs{x})^N},
\end{eqnarray*}
for all $N>d$.

\hfill $\blacksquare$ \vskip .5cm   

The definitions of $E_1,E_2,E_3$ in Lema \ref{l:conv_shrlts-wvlts}
allow us to have a ``height" of $2^{-3j}$ and a decreasing of
$(1+2^j\abs{x})^{-N}$.

\vskip .5cm The next proof regards the ``almost orthogonality" in
the Fourier domain.

\noindent\textbf{Proof of Lemma \ref{l:Overlap_bnd}}. This is a
direct consequence of the construction and dilation of the
shearlets. Suppose $d=2$. Since $k$ and $n$ are translation
parameters they do not seem reflected in the support of
$(\psi_{j,\ell,k})^\wedge$ or $(\psi_{i,m,n})^\wedge$. By
construction and by (\ref{e:Discrt_Shrlt_Cond_Cone_1}) one scale $j$
intersects with scales $j-1$ and $j+1$, only.

1) For one fixed scale $j$ and by (\ref{e:Discrt_Shrlt_Cond_Cone_2})
there exist \textbf{2} overlaps at the same scale $j$: one with
$(\psi_{j,\ell-1,k'})^\wedge$ and other with
$(\psi_{j,\ell+1,k''})^\wedge$ for all $k',k''\in\mathbb{Z}^2$.

2) Regarding scale $j-1$, one fixed $(\psi_{j,\ell,k})^\wedge$
overlaps with \textbf{3} other shearlets $(\psi_{j-1,m,k'})^\wedge$
at most for all $k,k'\in\mathbb{Z}^2$ because of 1) and because the
supports of the shearlets at scale $j-1$ have larger width than
those of scale $j$.

3) For a fixed scale $j$ consider the next three regions:
$\text{supp } (\psi_{j,\ell-1,k})^\wedge \cap \text{supp }
(\psi_{j,\ell,k'})^\wedge=R_{-1}$, $\text{supp }
(\psi_{j,\ell,k'})^\wedge \cap \text{supp }
(\psi_{j,\ell+1,k''})^\wedge=R_{+1}$ and $\text{supp
}(\psi_{j,\ell,k'})^\wedge \setminus (R_{-1}\cup R_{+1})=R_0$. Again
by construction, there can only be two overlaps for each $\xi$ at
any scale. Then, there exist at most two shearlets at scale $j+1$
that overlap with each of the three regions $R_i$, $i=-1,0,+1$ in
scale $j$: an aggregate of \textbf{6} for all translation parameters
$k,k',k''\in\mathbb{Z}^2$ at any scale $j$ or $j+1$.

\noindent Summing the number of overlaps at each scale gives the
result for $d=2$.

For the general case $d$ apply the same argument above for every
perpendicular direction to $\frak{d}$ .

\hfill $\blacksquare$ \vskip .5cm   

\textbf{Proof of Lemma \ref{l:Smplng_BndLim_Shrlt}.} Suppose first
that $g\in\mathcal{S}$. We can express $\hat{g}$ by its Fourier
series as
$$\hat{g}(\xi)=\sum_{k\in\mathbb{Z}^d} \abs{\text{det }A}^{-j/2} \mathbf{e}^{-2\pi i\xi A^{-j}B^{-\ell}k} \cdot
    \left(\int_{QB^\ell A^j} \hat{g}(\omega)\cdot\abs{\text{det }A}^{-j/2}\mathbf{e}^{2\pi i\omega A^{-j}B^{-\ell}k} d\omega\right).$$
By the Fourier inversion formula in $\hat{\mathbb{R}}^d$ we have
$$\hat{g}(\xi)=\sum_{k\in\mathbb{Z}^2} \abs{\text{det }A}^{-j/2} \mathbf{e}^{-2\pi i\xi A^{-j}B^{-\ell}k} \cdot
    g(A^{-j}B^{-\ell}k), \;\;\;\;\;\; \xi\in QB^\ell A^j.$$
Since $\hat{g}$ has compact support, $g(A^{-j}B^{-\ell}k)$ makes
sense (by the Paley-Wiener theorem). Since $\text{supp }
\hat{h}\subset QB^\ell A^j$ and $g\ast h=(\hat{g}\hat{h})^\vee$,
\begin{eqnarray*}
  g\ast h
    &=& \sum_{k\in \mathbb{Z}^d}\abs{\text{det }A}^{-j} g(A^{-j}B^{-\ell}k) [\mathbf{e}^{-2\pi i\xi A^{-j}B^{-\ell}k} \hat{h}(\cdot)]^\vee \\
    &=& \sum_{k\in \mathbb{Z}^d}\abs{\text{det }A}^{-j} g(A^{-j}B^{-\ell}k)
            h(x-A^{-j}B^{-\ell}k),
\end{eqnarray*}
which proves the convergence for $g\in\mathcal{S}$. To remove the
assumption $g\in\mathcal{S}$, we apply a standard regularization
argument to a $g\in\mathcal{S}'$ as done in p. 22 of \cite{Tr83} or
in Lemma A.4 of \cite{FJ90}. Let $\gamma\in\mathcal{S}$ satisfy
$\text{supp }\hat{\gamma}\subset B(0,1)$, $\hat{\gamma}(\xi)\geq 0$
and $\gamma(0)=1$. By Fourier inversion $\abs{\gamma(x)}\leq 1$ for
all $x\in\mathbb{R}^d$. For $0<\delta<1$, let
$g_\delta(x)=g(x)\gamma(\delta x)$. Then, $\text{supp
}\hat{g}_\delta$ is also compact, $g_\delta\in\mathcal{S}$,
$\abs{g_\delta}\leq \abs{g}$ for all $x\in\mathbb{R}^d$ and
$g_\delta\rightarrow g$ uniformly on compact sets as
$\delta\rightarrow 0$. Applying the previous result to $g_\delta$
and letting $\delta\rightarrow 0$ we obtain the result for general
$\mathcal{S}'$. This regularization argument (and, in fact, the
whole proof) is the same used in Lemma(6.10) in \cite{FJW}.

\hfill $\blacksquare$ \vskip .5cm   

\textbf{Proof of Lemma \ref{l:Plnchl-Polya_estimate_shrlts}.} By the
Paley-Wiener-Schwartz theorem $g$ is of exponential type, slowly
increasing and its point-wise values make sense (see Theorem 7.3.1
in \cite{Ho90}). Let $h\in\mathcal{S}$ satisfy supp
$\hat{h}\subseteq[-1,1]^d$ with $\hat{h}(\xi)=1$ for
$\xi\in[-\frac{1}{2},\frac{1}{2}]^d$. Write $g^y(x)=g(x+y)$. Then,
as in the proof of Lemma \ref{l:Smplng_BndLim_Shrlt}
$$g(x+y)=h_{A^{-j}B^{-\ell}}\ast g^y(x)=
    \sum_{\kappa\in\mathbb{Z}^d}\abs{\text{det }A}^{-j} g(A^{-j}B^{-\ell}\kappa+y)
    \abs{\text{det } A}^{j} h(B^\ell A^jx-\kappa).$$
Therefore, for any $y\in Q_{j,\ell,k}$,
$$\sup_{z\in Q_{j,\ell,k}} \abs{g(z)}\leq
    \sup_{\abs{x}<\text{diam} Q_{j,\ell,0}} \abs{g(x+y)}\leq
    \sum_{\kappa\in\mathbb{Z}^d} \abs{g(A^{-j}B^{-\ell}\kappa+y)}
    \sup_{\abs{x}<\text{diam} Q_{j,\ell,0}}\abs{h(B^\ell A^jx-\kappa)}.$$
But $h\in\mathcal{S}$ implies that for any $M>1$,
$$\sup_{\abs{x}<\text{diam} Q_{j,\ell,0}}\abs{h(B^\ell A^jx-\kappa)}
\leq \frac{C_M}{(1+\abs{\kappa})^M}.$$ Taking $M$ sufficiently large
and applying the $p$-triangle inequality $\abs{a+b}^p\leq
\abs{a}^p+\abs{b}^p$ if $0<p\leq 1$ or Hölder's inequality if $p>1$,
we obtain for any $y\in Q_{j,\ell,k}$, 
$$\sup_{z\in Q_{j,\ell,k}}\abs{g(z)}^p\leq C_p\sum_{\kappa\in\mathbb{Z}^d}
    \frac{\abs{g(A^{-j}B^{-\ell}\kappa+y)}^p}{(1+\abs{\kappa})^{d+1}}.$$
Integrating with respect to $y$ over the dyadic cube $Q_{j,k}$
yields
\begin{eqnarray*}
  2^{-jd}\sup_{z\in Q_{j,\ell,k}} \abs{g(z)}^p
    &\leq& C_p \sum_{\kappa\in \mathbb{Z}^d} (1+\abs{\kappa})^{-(d+1)}
        \int_{Q_{j,k}} \abs{g(A^{-j}B^{-\ell}\kappa+y)}^p dy
\end{eqnarray*}
Summing over $k\in\mathbb{Z}^d$,
\begin{eqnarray*}
  \abs{Q_{j}}^{\frac{d}{d+1}} \sum_{k\in\mathbb{Z}^d}\sup_{z\in Q_{j,\ell,k}} \abs{g(z)}^p
    &\leq& C_p \sum_{\kappa\in \mathbb{Z}^d} (1+\abs{\kappa})^{-(d+1)}
        \int_{\mathbb{R}^d} \abs{g(A^{-j}B^{-\ell}\kappa+y)}^p dx\\
    &=& C_p \norm{g}_{L^p}^p \sum_{\kappa\in \mathbb{Z}^d}
        (1+\abs{\kappa})^{-(d+1)}
    = C'_p \norm{g}_{L^p}^p,
\end{eqnarray*}
which finishes the proof.

\hfill $\blacksquare$ \vskip .5cm   

\vskip 0.5cm
\subsection{Proofs for Section \ref{S:AB-T-L}}\label{sS:Proofs_AB-T-L}
To prove our results we follow \cite{FJ90}, \cite[\S 6.3]{HW96},
\cite{BH05} and \cite[\S 1.3]{Tr83}. Some previous well known
definitions and results are necessary.

\begin{defi}\label{d:Max_Fcn}
For a function $g$ defined on $\mathbb{R}^d$ and for a real number
$\lambda>0$ the \textbf{Peetre's maximal function} (see Lemma 2.1 in
\cite{Pe75}) is
$$g^\ast_\lambda(x)=\sup_{y\in\mathbb{R}^d}\frac{\abs{g(x-y)}}{(1+\abs{y})^{d\lambda}}, \;\;\; x\in\mathbb{R}^d.$$
\end{defi}

\begin{lem}\label{l:Bnd_diff_Max_Fcn}
Let $g\in \mathcal{S}'(\mathbb{R}^d)$ be such that $\text{supp
}(\hat{g})\subseteq \{\xi\in\hat{\mathbb{R}}^d: \abs{\xi}\leq R\}$
for some $R>0$. Then, for any real $\lambda>0$ there exists a
$C_\lambda>0$ such that, for $\abs{\alpha}=1$,
$$(\partial^\alpha g)^\ast_\lambda(x)\leq C_\lambda g^\ast_\lambda(x), \;\;\; x\in\mathbb{R}^d.$$
\end{lem}   
\textbf{Proof.} Since $g\in\mathcal{S}'$ has compact support in the
Fourier domain, $g$ is regular. More precisely, by the
Paley-Wiener-Schwartz theorem $g$ is slowly increasing (at most
polinomialy) and infinitely differentiable (e.g., Theorem 7.3.1 in
\cite{Ho90}). Let $\gamma$ be a function in the Schwartz class such
that $\hat{\gamma}(\xi)=1$ if $\abs{\xi}\leq R$. Then,
$\hat{\gamma}(\xi) \hat{g}(\xi) =\hat{g}(\xi)$ for all
$\xi\in\hat{\mathbb{R}}^d$. Hence, $\gamma\ast g=g$ and
$\partial^\alpha g=\partial^\alpha\gamma\ast g$. Moreover,
\begin{eqnarray*}
  \abs{\partial^\alpha g(x-y)}
    &=& \abs{\int_{\mathbb{R}^d} \partial^\alpha \gamma(x-y-z)g(z)dz}
        = \abs{\int_{\mathbb{R}^d} \partial^\alpha\gamma(w-y)g(x-w)dw} \\
    &\leq& \int_{\mathbb{R}^d} \abs{\partial^\alpha\gamma(w-y)}
    (1+\abs{w-y})^{d\lambda}
        (1+\abs{y})^{d\lambda} \frac{\abs{g(x-w)}}{(1+\abs{w})^{d\lambda}}dw,
\end{eqnarray*}
because of the triangular inequality. Therefore,
$$\abs{\partial^\alpha g(x-y)}\leq
g^\ast_\lambda(x)(1+\abs{y})^{d\lambda}
     \int_{\mathbb{R}^d} \abs{\partial^\alpha\gamma(w-y)} (1+\abs{w-y})^{d\lambda} dw.$$
Since $\gamma\in\mathcal{S}$, the last integral equals a finite
constant $c_\lambda$, independent of $y$, and we obtain
$$\abs{\partial^\alpha g(x-y)}\leq c_\lambda g^\ast_\lambda(x)(1+\abs{y})^{d\lambda},$$
which shows the desired result.
\hfill $\blacksquare$ \vskip .5cm   

We have a relation between the Hardy-Littlewood maximal function
$\mathcal{M}(\abs{g}^{1/\lambda})(x)$ and the Peetre's maxmal
function $g^\ast_\lambda$.
\begin{lem}\label{l:Max_Fcn_leq_HL-Max-Fcn}
Let $\lambda>0$ and $g\in\mathcal{S}'$ be such that $\text{supp
}(\hat{g})\subseteq \{\xi\in\hat{\mathbb{R}}^d: \abs{\xi}\leq R\}$
for some $R>0$. Then, there exists a constant $C_\lambda>0$ such
that
$$g^\ast_\lambda (x)\leq C_\lambda \left(\mathcal{M}(\abs{g}^{1/\lambda})(x)\right)^\lambda, \;\;\; x\in\mathbb{R}^d.$$
\end{lem}
\textbf{Proof.} Since $g$ is band-limited, $g$ is differentiable on
$\mathbb{R}^d$ (by the Paley-Wiener-Schwartz theorem), so we can
consider the pointwise values of $g$. Let $x,y\in\mathbb{R}^d$ and
$0<\delta<1$. Choose $z\in\mathbb{R}^d$ such that $z\in
B_\delta(x-y)$. We apply the mean value theorem to $g$ and the
endpoints $x-y$ and $z$ to get
$$\abs{g(x-y)}\leq\abs{g(z)}+\delta\sup_{w:w\in B_\delta(x-y)}
( \abs{\nabla g(w)}).$$ Taking the $(1/\lambda)^{\text{th}}$ power
and integrating with respect to the variable $z$ over
$B_\delta(x-y)$, we obtain
\begin{eqnarray}\label{e:l:Max_Fcn_leq_HL-Max-Fcn}
\nonumber
  \abs{g(x-y)}^{1/\lambda}
    &\leq& \frac{c_\lambda}{\abs{B_\delta(x-y)}}
      \int_{B_\delta(x-y)} \abs{g(z)}^{1/\lambda}dz \\
    &&+ c_\lambda \delta^{1/\lambda} \sup_{w:w\in B_\delta(x-y)}
    (\abs{\nabla g(w)})^{1/\lambda}.
\end{eqnarray}
Since $B_\delta(x-y)\subset B_{\delta+\abs{y}}(x)$,
$$\int_{B_\delta(x-y)} \abs{g(z)}^{1/\lambda}dz \leq \int_{B_{\delta+\abs{y}}(x)}
\abs{g(z)}^{1/\lambda}dz
    \leq \abs{B_{\delta+\abs{y}}(x)} \mathcal{M}(\abs{g}^{1/\lambda})(x),$$
and the $\sup$ term on the right hand side of
(\ref{e:l:Max_Fcn_leq_HL-Max-Fcn}) is bounded by
\begin{eqnarray*}
  \sup_{w:w\in B_{\delta+\abs{y}}(x)} (\abs{\nabla g(w)})^{1/\lambda}
    &=& \sup_{t:\abs{t}<\delta+\abs{y}} (\abs{\nabla g(x-t)})^{1/\lambda} \\
    &\lesssim& (1+\delta+\abs{y})^{d} \left[(\nabla g)^\ast_\lambda(x)\right]^{1/\lambda}.
\end{eqnarray*}
Substituting these last two inequalities in
(\ref{e:l:Max_Fcn_leq_HL-Max-Fcn}) yields
\begin{eqnarray*}
  \abs{g(x-y)}^{1/\lambda}
    &\leq& c_\lambda\frac{\abs{B_{\delta+\abs{y}}(x)}}{\abs{B_\delta(x-y)}}\mathcal{M}(\abs{g}^{1/\lambda})(x)  \\
    && + c_\lambda\delta^{1/\lambda}(1+\delta+\abs{y})^d \left[(\nabla g)^\ast_\lambda(x)\right]^{1/\lambda},
\end{eqnarray*}
and since
$\abs{B_{\delta+\abs{y}}(x)}/\abs{B_\delta(x-y)}=(\delta+\abs{y})^d/\delta^d$,
we get
\begin{eqnarray*}
  \abs{g(x-y)}^{1/\lambda}
    &\leq& c_\lambda\frac{(\delta+\abs{y})^d}{\delta^d}\mathcal{M}(\abs{g}^{1/\lambda})(x)  \\
    && + c_\lambda\delta^{1/\lambda}(1+\delta+\abs{y})^d \left[(\nabla g)^\ast_\lambda(x)\right]^{1/\lambda}.
\end{eqnarray*}
Taking the $\lambda^{\text{th}}$ power yields
$$\frac{\abs{g(x-y)}}{(1+\abs{y})^{d\lambda}}\leq
c'_\lambda\left\{\frac{1}{\delta^{d\lambda}}[\mathcal{M}(\abs{g}^{1/\lambda})(x)]^\lambda
    + \delta\left[(\nabla g)^\ast_\lambda(x)\right] \right\},$$
since $\delta<1$ implies $(1+\delta+\abs{y})\leq 2(1+\abs{y})$.
Taking $\delta$ small enough so that $c'_\lambda
C_\lambda\delta<1/(2d)$ (where $C_\lambda$ is the constant in Lemma
\ref{l:Bnd_diff_Max_Fcn}) we obtain
$$g^\ast_\lambda(x)\leq
c_\lambda[\mathcal{M}(\abs{g}^{1/\lambda})(x)]^{\lambda}
    + \frac{1}{2}g^\ast_\lambda(x).$$
Assume for the moment that $g\in\mathcal{S}$, hence
$g^\ast_\lambda(x)<\infty$. So, we can subtract the second term in
the right-hand side of the previous inequality from the left-hand
side of the previous inequality and complete the proof for
$g\in\mathcal{S}$. To remove this assumption one uses the same
standard regularization argument as in the proof of Lemma
\ref{l:Smplng_BndLim_Shrlt}.
\hfill $\blacksquare$ \vskip .5cm   

\noindent Lemmata \ref{l:Bnd_diff_Max_Fcn} and
\ref{l:Max_Fcn_leq_HL-Max-Fcn} are Peetre's inequality for
$f\in\mathcal{S}'$ whose proofs can be found in the references above
and we reproduce them for completeness.

\vskip0.5cm \noindent\textbf{Proof of Lemma
\ref{l:SupConv_leq_HLMaxFcn}.} Let $g(x)=(\psi_{A^{-j}B^{-\ell}}\ast
f)(x)$. Since $\psi$ is band-limited, so is $g$. On one hand, $j\geq
0$ implies $C_d2^j\abs{y}\leq\abs{B^\ell A^j y}$, by Lemma
\ref{l:Nested_ellipsoids_HD}. Thus,
\begin{eqnarray*}
  g^\ast_\lambda(t)
    &=& \sup_{y\in\mathbb{R}^d}\frac{\abs{g(t-y)}}{(1+\abs{y})^{d\lambda}}
        \geq \sup_{y\in\mathbb{R}^d} \frac{\abs{(\psi_{A^{-j}B^{-\ell}}\ast
            f)(t-y)}}{(1+2^j\abs{y})^{d\lambda}} \\
    &=& \sup_{y\in\mathbb{R}^d} \frac{\abs{(\psi_{A^{-j}B^{-\ell}}\ast
            f)(t-y)}}{2^{d\lambda}(2^{-1}+2^{j-1}\abs{y})^{d\lambda}} \\
    &\geq& C_{d,\lambda} \sup_{y\in\mathbb{R}^d} \frac{\abs{(\psi_{A^{-j}B^{-\ell}}\ast
        f)(t-y)}}{(1+\abs{B^\ell A^jy})^{d\lambda}}
        = C_{d,\lambda} \abs{(\psi_{j,\ell,\lambda}^{\ast\ast})(t)}.
\end{eqnarray*}
On the other hand,
\begin{eqnarray*}
  \mathcal{M}(\abs{g}^{1/\lambda})(t)
    &=& \sup_{r>0} \frac{1}{\abs{B_r(t)}} \int_{B_r(t)} \abs{(\psi_{A^{-j}B^{-\ell}}\ast f)(y)}^{1/\lambda} dy \\
    &=& \mathcal{M}(\abs{(\psi_{A^{-j}B^{-\ell}}\ast
        f}^{1/\lambda})(t).
\end{eqnarray*}
The result follows from Lemma \ref{l:Max_Fcn_leq_HL-Max-Fcn} with
$t=x$.
\hfill $\blacksquare$ \vskip0.5cm   


To prove Lemma \ref{l:s_ast_bnd_s} we need the next result.
\begin{lem}\label{l:bnd_conv_seqN-seq3-2}
Let $i\geq j\geq 0$ and $0<a\leq r$. Also, let $Q$ and $P$ be
identified with $(i,m,n)$ and $(j,\ell,k)$, respectively. Then, for
all $N>(d+1)r/a$, any sequence $\{s_P\}_{P\in\mathcal{Q}^{j,\ell}}$
of complex numbers and any $x\in Q$,
$$(s_{r,N}^\ast)_Q:=\left(\sum_{P\in\mathcal{Q}^{j,\ell}}\frac{\abs{s_P}^r}{(1+2^j\abs{x_Q-x_P})^{N}}\right)^{1/r}
    \leq C_{a,r,d}
    \left[\mathcal{M}\left(\sum_{P\in\mathcal{Q}^{j,\ell}}\abs{s_P}^a\chi_P\right)(x)\right]^{1/a}.$$
Moreover, when $i=j$,
$$  \sum_{P\in\mathcal{Q}^{j,\ell}} \left[(s^\ast_{r,N})_P\tilde{\chi}_P(x)\right]^q
    \leq C_{a,r,d} \left[\mathcal{M}\left(\sum_{P\in\mathcal{Q}^{j,\ell}}
        (\abs{s_P}\tilde{\chi}_{P})^a\right)(x)\right]^{q/a}.$$
\end{lem}
\textbf{Proof.} Identify $(i,m,n)$ and $(j,\ell,k)$ with $Q$ and
$P$, respectively. Then, $x_Q=A^{-i}B^{-m}n$ and
$x_P=A^{-j}B^{-\ell}k$. Let
$\mathcal{Q}^{j,\ell}:=\{Q_{j,\ell,k}:k\in\mathbb{Z}^d\}$, then
$\mathcal{Q}^{j,\ell}$ is a partition of $\mathbb{R}^d$. Write
$l_P=\abs{x_Q-x_P}$. Thus, we bound the sum in the definition of
$(s_{r,N}^\ast)_Q$ as
$$\sum_{P\in\mathcal{Q}^{j,\ell}}\frac{\abs{s_P}^r}{(1+2^j\abs{x_Q-x_P})^N}
    \leq \left(\sum_{P\in\mathcal{Q}^{j,\ell}:l_P\leq 1} + \sum_{P\in\mathcal{Q}^{j,\ell}:l_P> 1}\right)
        \frac{\abs{s_P}^r}{(1+2^jl_P)^N}.$$
Choose $\lambda$ such that $N>(d+1)\lambda/d>(d+1)r/a$. Then, the
inequality $(2^j l_P)^N> (2^{(d+1)j/d}l_P)^\lambda$ holds whenever
$l_P>2^{j((d+1)\lambda/d-N)/(N-\lambda)}$. So, the previous
inequality is bounded by
$$\sum_{P\in\mathcal{Q}^{j,\ell}:l_P\leq 1} \frac{\abs{s_P}^r}{(1+2^jl_P)^N}
    + \sum_{P\in\mathcal{Q}^{j,\ell}:l_P> 1} \frac{\abs{s_P}^r}{(2^{(d+1)j/d}l_P)^\lambda}.$$
Defining
\begin{eqnarray*}
  D_0
    &=& \{k\in\mathbb{Z}^d:\abs{A^{-i}B^{-m}n-A^{-j}B^{-\ell}k}\leq 1\} \\
    &=& \{P\in\mathcal{Q}^{j,\ell}:l_P=\abs{x_Q-x_P}\leq 1\}
\end{eqnarray*}
and
\begin{eqnarray*}
  D_\nu
    &=& \{k\in\mathbb{Z}^d:2^{\nu-1}<2^{(d+1)j/d}\abs{A^{-i}B^{-m}n-A^{-j}B^{-\ell}k}\leq 2^\nu\} \\
    &=& \{P\in\mathcal{Q}^{j,\ell}:2^{\nu-1}<2^{(d+1)j/d}l_P\leq 2^\nu\},\;\;\;
    \nu=1,2,3,\ldots,
\end{eqnarray*}
we have that
\begin{eqnarray*}
  \sum_{P\in\mathcal{Q}^{j,\ell}}\frac{\abs{s_P}^r}{(1+2^jl_P)^N}
    &\leq& \sum_{P\in D_0}\abs{s_P}^r +2^\lambda\sum_{\nu=1}^\infty\sum_{P\in D_\nu}\frac{\abs{s_P}^r}{2^{\nu\lambda}} \\
    &\leq& 2^\lambda\sum_{\nu=0}^\infty \sum_{P\in D_\nu} \frac{\abs{s_P}^r}{2^{\nu\lambda}}
        \leq 2^\lambda\sum_{\nu=0}^\infty 2^{-\nu\lambda}
            \left(\sum_{P\in D_\nu}\abs{s_P}^a\right)^{r/a},
\end{eqnarray*}
because $1+2^jl_P\geq 1$ and $a\leq r$. Now, when $x\in Q_{i,m,n}=Q$
and $P\in D_\nu$ then, by the definition of $D_\nu$,
$P=Q_{j,\ell,k}\subset B_{(d+1)2^{\nu-(d+1)j/d}}(x)$
 (this holds because for $j\geq 0$ the
diameter of any $P=Q_{j,\ell,k}$ is less than $(d+1)$ and the
intervals in the definition of the $D_\nu$'s are dyadic with
$\nu\geq 0$). Thus,
$$\sum_{P\in D_\nu}\abs{s_P}^a=\abs{P}^{-1}\int_{B_{(d+1)2^{\nu-(d+1)j/d}}(x)}
    \sum_{P\in D_\nu}\abs{s_P}^a\chi_{P}(y)dy.$$
Hence, writing
$\abs{\tilde{B}}=\abs{B_{(d+1)2^{\nu-(d+1)j/d}}(x)}=C_d
(d+1)^d2^{d\nu-(d+1)j}$ we have that for $x\in Q_{i,m,n}=Q$,
\begin{eqnarray*}
  \sum_{P\in\mathcal{Q}^{j,\ell}}\frac{\abs{s_P}^r}{(1+2^jl_P)^N}
    &\leq& 2^\lambda\sum_{\nu=0}^\infty 2^{-\nu\lambda}\left( \frac{\abs{P}^{-1}\abs{\tilde{B}}}{\abs{\tilde{B}}}
        \int_{\tilde{B}}\sum_{P\in D_\nu}\abs{s_P}^a\chi_{P}(y)dy \right)^{r/a} \\
    &\leq& C_{a,r,d}\sum_{\nu=0}^\infty 2^{-\nu\lambda}2^{d\nu r/a} \left(\mathcal{M}
        \left(\sum_{P\in D_\nu}\abs{s_P}^a\chi_{P}\right)(x)\right)^{r/a}\\
    &\leq& C_{a,r,d}
        \left(\mathcal{M}\left(\sum_{P\in\mathcal{Q}^{j,\ell}}\abs{s_P}^a\chi_{P}\right)(x)\right)^{r/a},
\end{eqnarray*}
because $\abs{P}^{-1}=2^{(d+1)j}$ and $\lambda>dr/a$. 

To prove the second inequality multiply both sides by
$\tilde{\chi}_Q(x)$, rise to the power $q$ and sum over $Q\in
\mathcal{Q}^{j,\ell}$ to get
\begin{eqnarray*}
  \sum_{Q\in\mathcal{Q}^{j,\ell}} \left[(s^\ast_{r,N})_Q\tilde{\chi}_Q(x)\right]^q
    &\leq& C\sum_{Q\in\mathcal{Q}^{j,\ell}}  \left[\mathcal{M}\left(\sum_{P\in\mathcal{Q}^{j,\ell}}
        \abs{s_P}^a\chi_{P}\right)(x)\right]^{q/a}\tilde{\chi}_Q^q(x)\\
    &=& C\sum_{Q\in\mathcal{Q}^{j,\ell}}  \left[\mathcal{M}\left(\sum_{P\in\mathcal{Q}^{j,\ell}}
        (\abs{s_P}\tilde{\chi}_{P})^a\right)(x)\right]^{q/a}\chi_Q(x) \\
    &=& C \left[\mathcal{M}\left(\sum_{P\in\mathcal{Q}^{j,\ell}}
        (\abs{s_P}\tilde{\chi}_{P})^a\right)(x)\right]^{q/a},
\end{eqnarray*}
since $\mathcal{Q}^{j,\ell}$ is a partition of $\mathbb{R}^d$.

\hfill $\blacksquare$ \vskip .5cm   

\noindent\textbf{Proof of Lemma \ref{l:s_ast_bnd_s}}. Let $\lambda$
be such that $N>(d+1)\lambda/d>(d+1)\max(1,r/q,r/p)$. If
$r<\min(q,p)$, choose $a=r$. Otherwise, if $r\geq\min(q,p)$, choose
$a$ such that $r/(\lambda/d)<a<\min(r,q,p)$. It is always possible
to choose such an $a$ since $\lambda/d>\max(1,r/q,r/p)$ implies
$r/(\lambda/d)<\min(r,q,p)$. In both cases we have that
$$0<a\leq r<\infty, \;\;\; \lambda>dr/a, \;\;\; q/a>1, \;\;\; p/a>1.$$
The previous argument is similar to that in \cite{BH05}. Then, by
Lemma \ref{l:bnd_conv_seqN-seq3-2} and Theorem \ref{t:Feff-Stn_Ineq}
\begin{eqnarray*}
  \norm{\mathbf{s}^\ast_{r,N}}_{\mathbf{f}^{\alpha,q}_p(AB)}
    &=& \norm{\left(\sum_{P\in\mathcal{Q}_{AB}}(\abs{P}^{-\alpha}(s^\ast_{r,N})_P\tilde{\chi}_P)^q\right)^{1/q}}_{L^p} \\
    &=& \norm{\left(\sum_{j\geq 0}\sum_{\abs{[\ell]}\preceq2^j}\abs{Q_{j,\ell}}^{-\alpha q}\sum_{P\in\mathcal{Q}^{j,\ell}}
        (s^\ast_{r,N})_P^q\tilde{\chi}_P^q \right)^{1/q}}_{L^p} \\
    &\leq& C_{a,r,d}\norm{\left(\sum_{j\geq 0}\sum_{\abs{[\ell]}\preceq2^j}\abs{Q_{j,\ell}}^{-\alpha q}
        \left[\mathcal{M}\left(\sum_{P\in \mathcal{Q}^{j,\ell}}(\abs{s_P}\tilde{\chi}_{P})^a\right)\right]^{q/a} \right)^{1/q}}_{L^p} \\
    &=& C_{a,r,d}\norm{\left(\sum_{j\geq 0}\sum_{\abs{[\ell]}\preceq2^j}
        \left[\mathcal{M}\left(\sum_{P\in \mathcal{Q}^{j,\ell}}(\abs{P}^{-\alpha}\abs{s_P}\tilde{\chi}_{P})^a\right)\right]^{q/a} \right)^{1/q}}_{L^p} \\
    &=& C_{a,r,d}\norm{\left(\sum_{j\geq 0}\sum_{\abs{[\ell]}\preceq2^j}
        \left[\mathcal{M}\left(\sum_{P\in \mathcal{Q}^{j,\ell}}(\abs{P}^{-\alpha}\abs{s_P}\tilde{\chi}_{P})^a\right)\right]^{q/a}
        \right)^{a/q}}_{L^{p/a}}^{1/a}\\
    &\leq& C_{a,r,d}\norm{\left(\sum_{j\geq 0}\sum_{\abs{[\ell]}\preceq2^j} \left(\sum_{P\in \mathcal{Q}^{j,\ell}}
        (\abs{P}^{-\alpha}\abs{s_P}\tilde{\chi}_P )^a\right)^{q/a}\right)^{a/q}}_{L^{p/a}}^{1/a} \\
    &=& C_{a,r,d}\norm{\left(\sum_{j\geq 0}\sum_{\abs{[\ell]}\preceq2^j}\sum_{P\in \mathcal{Q}^{j,\ell}}
        (\abs{P}^{-\alpha}\abs{s_P}\tilde{\chi}_P )^q\right)^{a/q}}_{L^{p/a}}^{1/a} \\
    &=& C_{a,r,d}\norm{\left(\sum_{j\geq 0}\sum_{\abs{[\ell]}\preceq2^j} \sum_{P\in \mathcal{Q}^{j,\ell}}
        (\abs{P}^{-\alpha}\abs{s_P}\tilde{\chi}_P )^q\right)^{1/q}}_{L^{p}} \\
    &=& C_{a,r,d}\norm{\mathbf{s}}_{f^{\alpha,q}_p(AB)},
\end{eqnarray*}
because $\mathcal{Q}^{j,\ell}$ is a partition of $\mathbb{R}^d$.

The reverse inequality is trivial since $\abs{s_Q}\leq
(s^\ast_{r,N})_Q$ always holds.
\hfill $\blacksquare$ \vskip .5cm   

%

%


\end{document}